\documentclass[11pt,a4paper,reqno]{amsart}
\usepackage[margin=1in]{geometry}
\usepackage{amsmath,amssymb,amsthm,mathtools}
\usepackage[T1]{fontenc}
\usepackage{lmodern,microtype}
\usepackage[colorlinks=true,linkcolor=blue,citecolor=red,urlcolor=blue]{hyperref}

\allowdisplaybreaks[2]

\newtheorem{theorem}{Theorem}[section]
\newtheorem{lemma}[theorem]{Lemma}
\newtheorem{proposition}[theorem]{Proposition}
\newtheorem{corollary}[theorem]{Corollary}
\theoremstyle{remark}
\newtheorem{remark}[theorem]{Remark}

\newcommand{\R}{\mathbb R}
\newcommand{\dd}{\,\mathrm d}
\newcommand{\la}[1]{\langle #1\rangle}
\newcommand{\tr}{\operatorname{tr}}
\newcommand{\diver}{\operatorname{div}}

\newcommand{\diag}{\operatorname{diag}}
\newcommand{\cD}{\mathcal D}
\newcommand{\cE}{\mathcal E}

\newcommand{\cG}{\mathcal G}
\newcommand{\cR}{\mathcal R}
\newcommand{\cK}{\mathcal K}

\numberwithin{equation}{section}
\title[Macroscopic $L^\infty$ bounds for the Landau--Coulomb equation]{Macroscopic $L^\infty$ bounds for the homogeneous\\
Landau--Coulomb equation}

\author{Lingda Xu}
\address{Zhizhen School of Interdisciplinary Mathematical Sciences, The Chinese University of Hong Kong, Shatin, N.T., Hong Kong SAR, China}
\email{xuldmath@gmail.com}


\hypersetup{
  pdftitle={Macroscopic L-infinity bounds for the homogeneous Landau--Coulomb equation},
  pdfauthor={Lingda Xu}
}

\begin{document}

\begin{abstract}
We establish an instantaneous global $L^\infty(\R^3)$ bound at positive times for smooth, rapidly decaying solutions to the three-dimensional spatially homogeneous Landau equation with Coulomb potential, where the upper bound depends solely on the initial mass, momentum, energy, and entropy. This resolves an open problem posed by Golding, Gualdani, and Loher \cite{GGL} without requiring any initial higher velocity moments, Fisher information, or Lebesgue integrability. The proof builds on two main analytical ingredients. First, we establish the instantaneous generation of the negatively weighted Fisher information with weight $\la{v}^{-1}$. On the lifted two-particle phase space, we compare the separable weight matrix against a collision-invariant pair weight and control the rotational commutator near the diagonal via the quadratic local mass decay supplied by the conserved energy, closing a Riccati differential inequality at the second-moment threshold. Second, to overcome the $\la{v}^{-3}$ radial degeneracy of the diffusion matrix at large velocities, we prove a truncated anisotropic mixed Nash inequality on upper level sets of the density that recovers longitudinal regularity from two transverse Fisher directions and the nonlinear high-density part of the Coulomb entropy dissipation. Localizing this functional inequality to radially aligned velocity cells with the critical weight $\la{v}^{3/2}$ and coupling it with a weighted truncated entropy identity across nested exterior domains yields a superlinear De Giorgi recurrence that bounds the velocity tail uniformly in $\R^3$.
\end{abstract}

\maketitle

\section{Introduction}

\subsection{Problem formulation and main results}

The spatially homogeneous Landau equation models collisional relaxation of a plasma in velocity space through grazing binary interactions. In three dimensions, the Coulomb potential represents the most singular physical interaction: the collision kernel diverges at small relative velocities, while the induced diffusion matrix degenerates anisotropically at large velocities. We consider nonnegative, smooth, rapidly decaying solutions $f=f(t,v)$ on $[0,T]\times\R^3$ to the Cauchy problem
\begin{equation}\label{eq:landau}
 \partial_t f=Q(f,f)=\nabla_v\cdot\bigl(A[f]\nabla_v f-f\nabla_v a[f]\bigr),
\end{equation}
where the diffusion matrix $A[f]$ and the scalar potential $a[f]$ are given by the convolutions
\begin{equation}\label{eq:kernels}
 \begin{aligned}
 A[f](v)&=\int_{\R^3}K(v-w)f(w)\dd w,
 &K(z)&=\frac{1}{8\pi|z|}\bigl(I_3-\tilde z\otimes\tilde z\bigr),\\
 a[f](v)&=\frac{1}{4\pi}\int_{\R^3}\frac{f(w)}{|v-w|}\dd w,
 &\tilde z&=\frac{z}{|z|}\quad (z\ne 0).
 \end{aligned}
\end{equation}
Here $I_3$ is the $3\times 3$ identity matrix. Differentiating the kernels in \eqref{eq:kernels} in the sense of distributions yields the structural relations
\[
 \tr A[f]=a[f],\qquad \diver A[f]=\nabla a[f],\qquad -\Delta a[f]=f,
\]
which allow \eqref{eq:landau} to be written in the non-divergence form
\[
 \partial_t f=A[f]:D^2 f+f^2.
\]
By translation and scaling invariance, we assume throughout the paper that the initial mass, momentum, and energy are normalized and that the initial Boltzmann entropy is bounded above by $H$:
\begin{equation}\label{eq:normalization}
 \int_{\R^3} f(t,v)\dd v=1,\quad \int_{\R^3} vf(t,v)\dd v=0,\quad
 \int_{\R^3} |v|^2f(t,v)\dd v=3,\quad
 \int_{\R^3} f(0,v)\log f(0,v)\dd v\le H.
\end{equation}
Along the flow of \eqref{eq:landau}, the mass, momentum, and energy are conserved, while the Boltzmann entropy is nonincreasing in time. Constants denoted by $C_H$ and $c_H$ depend only on the entropy upper bound $H$ and the normalization \eqref{eq:normalization}. Throughout the paper, we write $A\sim B$ to indicate that $c^{-1}B\le A\le cB$ for a universal constant $c>1$.

Writing $\la{v}=(1+|v|^2)^{1/2}$, we define the weighted Fisher information and the weighted logarithmic Hessian functional by
\begin{equation}\label{eq:IJ}
 I_{-q}(f)=\int_{\R^3} f\la{v}^{-q}|\nabla\log f|^2\dd v,\qquad
 J_4(f)=\int_{\R^3} f\la{v}^{-4}|D^2\log f|^2\dd v,
\end{equation}
where matrix norms denote Frobenius norms, the integrand $f|\nabla\log f|^2$ is interpreted as $4|\nabla\sqrt{f}|^2$, and $I(f)=I_0(f)$ is the standard unweighted Fisher information. A solution $f$ on $[0,T]\times\R^3$ is called smooth and rapidly decaying if, for every integer $N\ge0$ and every multi-index $\alpha$,
\[
 \sup_{0\le t\le T}\sup_{v\in\R^3}\la{v}^N
 |\partial_v^\alpha f(t,v)|<\infty.
\]
While these qualitative Schwartz seminorms are used to justify intermediate integrations by parts, all quantitative constants in our main estimates depend only on the macroscopic quantities in \eqref{eq:normalization}.

Our main theorem shows that any smooth, rapidly decaying solution to the Landau--Coulomb equation is instantaneously bounded in $L^\infty(\R^3)$ for $t>0$ by a constant depending only on the conserved quantities and the initial entropy.

\begin{theorem}[Global $L^\infty$ bound for positive times]\label{thm:global-bound}
For every $H\in\R$ and $0<\tau<T$, there exists a constant $C(H,\tau,T)>0$, depending only on $H$, $\tau$, and $T$, such that every smooth, rapidly decaying solution $f$ of \eqref{eq:landau}--\eqref{eq:normalization} satisfies
\[
 \sup_{\tau\le t\le T}\|f(t)\|_{L^\infty(\R^3)}
 \le C(H,\tau,T).
\]
\end{theorem}

Theorem~\ref{thm:global-bound} provides a quantitative a priori $L^\infty(\R^3)$ regularization estimate solely in terms of the macroscopic conservation laws and entropy, answering the open question raised by Golding, Gualdani, and Loher \cite{GGL}. Upgrading this a priori estimate to a global existence theorem for general initial data with finite mass, energy, and entropy requires an additional approximation and initial-trace stability argument, which we do not pursue here.

The first pillar in the proof of Theorem~\ref{thm:global-bound} is the instantaneous generation of the negatively weighted Fisher information $I_{-1}(f)$ at the second-moment threshold.

\begin{theorem}\label{thm:critical-fisher}
Under the normalization \eqref{eq:normalization}, there exist positive constants $c_H,C_H,\tilde c_H,\tilde C_H$ depending only on $H$ such that, for almost every $t\in(0,T)$,
\begin{equation}\label{eq:main-diff}
 \frac{\dd}{\dd t}I_{-1}(f)
 \le -c_HJ_4(f)+C_H\bigl(1+I_{-1}(f)\bigr)^{3/2}
 \le -\tilde c_HI_{-1}(f)^2+\tilde C_H.
\end{equation}
Consequently, for all $0<t\le T$,
\begin{equation}\label{eq:main-bound}
 I_{-1}(f(t))\le C_H\bigl(1+t^{-1}\bigr).
\end{equation}
\end{theorem}

Combining Theorem~\ref{thm:critical-fisher} with Sobolev embedding and potential estimates yields global $L^{9/5}(\R^3)$ integrability for the density, uniform bounds on the collision coefficients, and local $L^\infty$ bounds on bounded velocity balls.

\begin{corollary}\label{cor:smoothing}
Every solution as in Theorem~\ref{thm:critical-fisher} satisfies, for all $0<t\le T$,
\begin{equation}\label{eq:cor-bounds}
 \begin{aligned}
 \|f(t)\|_{L^{9/5}(\R^3)}&\le C_H\bigl(1+t^{-2/3}\bigr),\\
 \|a[f(t)]\|_{L^\infty(\R^3)}&\le C_H\bigl(1+t^{-1/2}\bigr),\\
 \|\nabla a[f(t)]\|_{L^{9/2}(\R^3)}&\le C_H\bigl(1+t^{-2/3}\bigr).
 \end{aligned}
\end{equation}
Furthermore, for every $0<\tau<T$ and every finite radius $R>0$,
\begin{equation}\label{eq:local-infty}
 \sup_{\tau\le t\le T}\|f(t)\|_{L^\infty(B_R)}
 \le C(H,\tau,T,R).
\end{equation}
\end{corollary}

While Corollary~\ref{cor:smoothing} controls the solution on any bounded velocity ball $B_R$, the local Moser iteration constant degenerates rapidly as $R\to\infty$. Sections~\ref{sec:truncated-mixed-nash}--\ref{sec:original-tail-DG} remove this radius dependence by establishing a uniform bound outside a fixed ball $B_{2R}$ of radius $R=R(H,\tau,T)$.

\subsection{Background and related works}

The global regularity theory for the spatially homogeneous Landau--Coulomb equation advanced decisively with the work of Guillen and Silvestre \cite{GS}. By lifting the collision operator to the two-particle phase space $\R^6$ and exploiting the rotational geometry of binary collisions, they proved that the unweighted Fisher information $I(f)$ is nonincreasing in time. This monotonicity rules out finite-time blowup and yields global smooth solutions whenever the initial data possess finite Fisher information and suitable spatial decay.

Following \cite{GS}, substantial effort has been devoted to relaxing the regularity and moment requirements on the initial data. Golding, Gualdani, and Loher \cite{GGL} established $L^\infty$ regularization estimates down to the critical weighted $L^{3/2}$ threshold. In the introduction of \cite{GGL}, they formulated the open problem of whether smooth solutions to the homogeneous Landau--Coulomb equation admit $L^\infty(\R^3)$ bounds for positive times depending solely on the physical mass, energy, and entropy, and separately whether the initial smoothing rate $t^{-1}$ is sharp. Theorem~\ref{thm:global-bound} resolves the first problem in the affirmative.

A closely related line of investigation concerns the instantaneous generation of Fisher information at $t>0$. For the unweighted functional $I(f)$, Ji \cite{JiDissipation} proved a dissipation bound controlling a logarithmic Hessian integral weighted by $\la{v}^{-5}$ in the Coulomb case, which closes a differential inequality for $I(f)$ under a fifth moment assumption. Independently, Desvillettes, Golding, Gualdani, and Loher \cite{DGGL} established the instantaneous production of $I(f)$ assuming a finite eighteenth moment. Weighted Fisher functionals in kinetic equations with soft potentials were earlier studied by Wu \cite{Wu}, and general polynomially weighted Fisher identities together with weighted Hessian interpolations were systematically developed by Feng and Wang \cite[Lemmas 3.7 and 3.14]{FW} in the context of Kac's program. Rather than seeking unweighted Fisher production via higher moments, Theorem~\ref{thm:critical-fisher} pairs the negatively weighted functional $I_{-1}(f)$ with the Hessian integral $J_4(f)$ so that the closing interpolation requires only the conserved second moment. Our analysis also builds on foundational results on entropy dissipation \cite{Desvillettes}, conditional parabolic regularization \cite{Silvestre}, uniqueness of bounded solutions \cite{Fournier}, and truncated entropy techniques for soft potentials \cite{GIJV}.

\subsection{Main difficulties and proof strategy}

Establishing a global $L^\infty(\R^3)$ bound under only the conservation laws and entropy requires overcoming two fundamental obstacles created by the grazing singularity at small relative velocities and the anisotropic $\la{v}^{-3}$ degeneracy at large velocities when no moments above order two are available.

\medskip
\noindent\textbf{Difficulties:}
\begin{enumerate}
\item \emph{Closing the Fisher dissipation inequality at the second-moment threshold.}
In the lifted framework of Guillen and Silvestre \cite{GS}, the time derivative of the Fisher information is expressed through the action of the spherical diffusion operator $\sum_{k=1}^3 X_k^2$ on the product density $F(v,w)=f(v)f(w)$. Here the divergence-free vector fields $X_k$ generate rotations in the relative velocity $z=v-w$ on each collision sphere of fixed center of mass $c=(v+w)/2$ and radius $r=|z|$. For the unweighted functional $I(f)$, Ji's sharp $\Gamma_2$ inequality on the real projective sphere \cite{JiSphere} absorbs the spherical curvature term directly. However, extracting coercivity from the spherical dissipation produces a logarithmic Hessian integral carrying a negative velocity weight. In our setting, the natural coercive Hessian functional is $J_4(f)$ with weight $\la{v}^{-4}$. Attempting to bound the unweighted Fisher information $I(f)$ by $J_4(f)$ via integration by parts gives
\[
 I(f)=-\int_{\R^3}f\Delta\log f\dd v
 \le \left(3J_4(f)\int_{\R^3}\la{v}^4 f\dd v\right)^{1/2},
\]
which requires a fourth velocity moment (and similarly a fifth moment in \cite{JiDissipation}, where the Hessian weight is $\la{v}^{-5}$).

To close the interpolation using only the conserved second moment, one is forced to replace $I(f)$ by the negatively weighted Fisher functional $I_{-1}(f)$ with weight $\omega(v)=\la{v}^{-1}$. On the lifted phase space $\R^6$, however, the corresponding separable weight matrix
\[
 S(v,w)=\diag\bigl(\la{v}^{-1}I_3,\la{w}^{-1}I_3\bigr)
\]
is no longer constant along the collision spheres. Differentiating $S(v,w)$ twice along the rotation fields $X_k$ destroys the exact spherical structure and generates a second-order commutator kernel with an $r^{-3}$ singularity near the collision diagonal $\{v=w\}$, coupling the scores and logarithmic Hessians of both particles.

\item \emph{Preventing polynomial weight loss in the degenerate velocity tail.}
Once Theorem~\ref{thm:critical-fisher} provides the bound on $I_{-1}(f(t))$ for $t>0$, three-dimensional Sobolev embedding applied to $\sqrt{\omega f}$ yields $\la{v}^{-1}f\in L^3(\R^3)$, which interpolates with the conserved energy to give $f\in L^{9/5}(\R^3)$ and $a[f]\in L^\infty(\R^3)$. On any bounded velocity ball $B_R$, the diffusion matrix $A[f]$ is uniformly elliptic and classical Moser iteration yields the local bound \eqref{eq:local-infty}. At large velocities, however, the lower ellipticity of $A[f](v)$ degenerates anisotropically at different rates in the radial and transverse directions:
\[
 A[f](v)\ge c_H\bigl(\la{v}^{-3}\Pi_v+\la{v}^{-1}\Pi_v^\perp\bigr),
 \qquad \Pi_v=\frac{v\otimes v}{|v|^2},\quad \Pi_v^\perp=I_3-\Pi_v.
\]
In standard parabolic De Giorgi or Moser iterations that linearize \eqref{eq:landau} into a scalar drift-diffusion equation (such as \cite{Silvestre,GGL}), one must rely on the weakest eigenvalue $\la{v}^{-3}$ in the radial direction to apply the three-dimensional Sobolev inequality. Each step of the iteration then loses a polynomial power of $\la{v}$. Over infinitely many steps on the exterior domain $\{|v|>R\}$, or even in a single step without higher moments, this radial weight loss requires velocity moments well above order two to close.
\end{enumerate}

\medskip
\noindent\textbf{Strategies and main innovations:}
\begin{enumerate}
\item \emph{Comparison with a collision-invariant pair weight and local mass decay near the diagonal (Sections~\ref{sec:lifted}--\ref{sec:coercivity}).}
To resolve the first difficulty, we decompose the separable weight matrix $S(v,w)$ around the pair weight
\[
 \Omega(v,w)=\left(1+\frac{|v|^2+|w|^2}{2}\right)^{-1/2}
 =\left(\la{c}^2+\frac{r^2}{4}\right)^{-1/2}.
\]
Because $|v|^2+|w|^2=2|c|^2+\frac{1}{2}r^2$ is conserved by binary collisions, $\Omega$ is constant on every collision sphere and is annihilated by all three rotation fields ($X_k\Omega=0$). Furthermore, the elementary inequality $\la{v}^2\le\la{v}^2+\la{w}^2$ guarantees the global lower bound
\[
 \frac{\omega(v)}{\Omega(v,w)}\ge\frac{1}{\sqrt{2}}\qquad\text{on }\R^6.
\]
Since Ji's projective sphere inequality \cite{JiSphere} provides an $11$-fold coercivity ratio between the spherical dissipation $D_\Omega$ and the curvature integral $R_\Omega$, and $\sqrt{2}>9/11$, replacing $S$ by $2^{-1/2}\Omega I_6$ in the leading dissipation preserves a strictly positive coercive margin:
\[
 \kappa_0=\sqrt{2}-\frac{9}{11}>0.
\]

It then remains to control the commutator generated by the perturbation matrix $T=S-\Omega I_6$, which vanishes identically on the collision diagonal $\{v=w\}$. We split $\R^6$ using a radial cutoff at the relative scale $r\sim\delta\la{c}$ for a small parameter $\delta>0$. For separated pairs with $r\ge\delta\la{c}$, the relative velocity $r$ dominates both $\la{v}$ and $\la{w}$, so the commutator kernel decays like $r^{-1}$ and is controlled directly by $C_\delta I_{-1}(f)$. For nearby pairs with $r\le 2\delta\la{c}$, expanding the second commutator directly would leave a singular $r^{-2}\Omega/\la{c}$ kernel. Instead, we reverse one rotational integration by parts and integrate by parts once in Euclidean coordinates. This transfers a derivative onto the logarithmic Hessian $D^2\log f(v)$ while gaining a linear Taylor factor of $r\le 2\delta\la{c}$ from the vanishing of $T$ on the diagonal. Because $\la{w}\sim\la{v}$ for $r\le 2\delta\la{c}$, Chebyshev's inequality and the conserved energy yield quadratic local mass decay:
\[
 \int_{\{r\le 2\delta\la{c}\}}f(w)\dd w\le C\la{v}^{-2},
 \qquad
 \int_{\{r\le 2\delta\la{c}\}}rf(w)\dd w\le C\delta\la{v}^{-1}.
\]
Combined with the kernel factor $\Omega/\la{c}^2\sim\la{v}^{-3}$, this local mass decay converts the Hessian error into $C\delta J_4(f)$, which is absorbed by $D_\Omega$ for small $\delta$ without invoking any higher moment. An exact single-particle integration by parts against $\omega(v)=\la{v}^{-1}$ then yields $I_{-1}(f)^2\le 12 J_4(f)$ and closes the Riccati inequality \eqref{eq:main-diff}.

\item \emph{Truncated anisotropic mixed Nash inequality and radially aligned tail De Giorgi iteration (Sections~\ref{sec:truncated-mixed-nash}--\ref{sec:original-tail-DG}).}
To overcome the $\la{v}^{-3}$ radial degeneracy without higher moments, we abandon the linear parabolic perspective in the velocity tail and retain the nonlinear high-density part of the Coulomb entropy dissipation. For any density threshold $K>0$ and truncated amplitude $h=(\sqrt{f}-\sqrt{K})_+^2$, we prove in Section~\ref{sec:truncated-mixed-nash} the homogeneous truncated mixed Nash inequality
\[
 \int_{\R^3}h^{5/3}\dd v
 \le C\left(\mu^{2/3}X+\mu^{1/2}X^{1/2}\cD^{1/3}\right),
\]
where the truncated mass $\mu$ and the two-dimensional transverse Fisher energy $X$ are given by
\[
 \mu=\int_{\{f>K\}}f\dd v,\qquad
 X=4\int_{\R^3}|\nabla_\perp\sqrt{h}|^2\dd v,
\]
and $\cD$ is the double-integral Coulomb entropy dissipation restricted to the high-density interaction set $\{f>K\}\times\{f>K\}$. The underlying geometric mechanism is that the angular projection $I_3-\tilde z\otimes\tilde z$ in the Landau kernel $K(v-w)$ couples the missing longitudinal derivative $\partial_{v_1}\log f$ with the transverse gradient $\nabla_\perp\log f$ whenever $v-w$ is not parallel to $e_1$. By introducing the generalized score $\sigma=\nabla(f-K)_+/(f\mathbf{1}_{\{f>K\}})$ and testing the cross product $(v-w)\times(\sigma(v)-\sigma(w))$ against the transverse moment $(w_j-\bar v_j)f(w)\mathbf{1}_{\{f>K\}}$, we recover longitudinal regularity on compact cells directly from $\cD$ and $X$ with zero boundary terms on the level set $\{f=K\}$.

In Section~\ref{sec:weighted-truncated-nash}, we localize this functional inequality to the exterior velocity domain by covering $\{|v|>R\}$ with unit-scale cells $\{\chi_j\}$ whose local longitudinal axis is aligned with the radial direction $v_j/|v_j|$ of each cell center. In this radially aligned frame, the two transverse directions experience the much stronger $\la{v}^{-1}$ background ellipticity of the low-density diffusion matrix, completely bypassing the degenerate $\la{v}^{-3}$ radial eigenvalue. Introducing the polynomial weight $W(v)=\la{v}^{3/2}$, we find that the critical exponent $\alpha=3/2$ exactly cancels the $\la{v_j}^{-1}$ transverse metric factor in the cell Sobolev embedding, yielding a localized weighted Nash inequality in the velocity tail with constants independent of the cell radius.

Finally, in Section~\ref{sec:original-tail-DG}, we combine this weighted Nash inequality with an exact evolution identity for the weighted truncated entropy
\[
 E_{p,k}(f)=\int_{\R^3}p(v)^2\bigl(f\log(f/k)-f+k\bigr)_+\dd v
\]
across nested exterior domains $\{|v|>R_n\}$ and thresholds $k_n=K_0(2-2^{-n})$. Exploiting the distributional identity $\diver\diver K=-\delta_0$, we show that integrating the high-density cross term by parts produces a diagonal delta mass $\int_{\R^3} p^2(f-k)_+^2\dd v$ that cancels the quadratic reaction term $+f^2$ in the Landau equation up to a linear excess term and a nonlocal weight commutator controlled solely by the weighted upper-level mass. Because the critical weight $W(v)=\la{v}^{3/2}$ lies strictly below the conserved second moment $\la{v}^2$, the bound $\|\la{v}^{-1}f\|_{L^3}<\infty$ from Theorem~\ref{thm:critical-fisher} provides small initial tail energy $Y_0,Y_1\le C_{H,\tau,T}R^{-1/4}$ at $t=\tau/4$, allowing a two-step superlinear De Giorgi recurrence to drive the tail entropy to zero at level $2K_0$.
\end{enumerate}

The remainder of the paper is organized as follows. Sections~\ref{sec:lifted}--\ref{sec:coercivity} establish the lifted matrix Fisher identity, the commutator estimates for nearby and separated pairs, and the coercivity of the pair dissipation, completing the proof of Theorem~\ref{thm:critical-fisher}. Section~\ref{sec:local} derives the global $L^{9/5}$ integrability, the collision coefficient bounds, and the local Moser iteration of Corollary~\ref{cor:smoothing}. Sections~\ref{sec:truncated-mixed-nash} and \ref{sec:weighted-truncated-nash} prove the unweighted and radially aligned weighted truncated mixed Nash inequalities. Finally, Section~\ref{sec:original-tail-DG} establishes the weighted truncated entropy identity and carries out the exterior De Giorgi iteration to complete the proof of Theorem~\ref{thm:global-bound}.

\section{The lifted Fisher identity}\label{sec:lifted}

We begin by deriving the lifted evolution identity for the weighted Fisher information on the two-particle phase space $\R^6=\R^3_v\times\R^3_w$. The differential identities in this section and Section~\ref{sec:commutator} are first computed for smooth, strictly positive densities with rapid decay; the approximation and spatial cutoff justification for general smooth, rapidly decaying solutions is given in Section~\ref{app:justification}.

Given a density $f(v)$ on $\R^3$, we introduce the product density $F(v,w)$, the single-particle and two-particle score vectors $s_v$ and $\tilde s$, the center-of-mass and relative-velocity coordinates $(c,z)$, the separable weight matrix $S(v,w)$, and the pair weight $\Omega(v,w)$:
\begin{equation}\label{eq:lift}
 \begin{aligned}
 F(v,w)&=f(v)f(w),&\tilde s(v,w)&=\nabla_{v,w}\log F(v,w)=(s_v,s_w),\\
 z&=v-w,\quad r=|z|,&c&=\frac{v+w}{2},\quad \rho=\la{c},\\
 \omega(v)&=\la{v}^{-1},&S(v,w)&=\diag\bigl(\omega(v)I_3,\omega(w)I_3\bigr),\\
 \Omega(v,w)&=\left(1+\frac{|v|^2+|w|^2}{2}\right)^{-1/2}
 =\left(\rho^2+\frac{r^2}{4}\right)^{-1/2}.\span\omit
 \end{aligned}
\end{equation}
For each coordinate index $k\in\{1,2,3\}$, let $e_k$ denote the $k$-th standard basis vector of $\R^3$. We define the lifted rotation vector fields $X_k$ and their associated score gradients by
\begin{equation}\label{eq:fields}
 \begin{aligned}
 b_k(z)&=e_k\times z,&\Xi_k(v,w)&=r^{-3/2}\bigl(b_k(z),-b_k(z)\bigr),&X_k&=\Xi_k\cdot\nabla_{v,w},\\
 x_k&=X_k\log F,&V_k&=\nabla_{v,w}x_k,&(L_k)_{ij}&=\partial_j(\Xi_k)_i.
 \end{aligned}
\end{equation}
Each vector field $X_k$ is divergence-free on $\R^6$ and generates rotations in the relative velocity $z$ while holding $c$ and $r$ fixed. Consequently, $X_k$ annihilates every function of $(c,r)$, and in particular $X_k\Omega=0$.

To separate the invariant pair weight $\Omega I_6$ from the perturbation $S-\Omega I_6$ in a linear fashion, we compute the lifted Fisher evolution associated with an arbitrary smooth symmetric $6\times 6$ matrix weight $M(v,w)$. For each $k\in\{1,2,3\}$, define the first and second Lie commutator matrices of $M$ along $X_k$ by
\begin{equation}\label{eq:RK}
 \cR_{M,k}=L_kM+ML_k^\top-X_kM,\qquad
 \cK_{M,k}=L_k\cR_{M,k}+\cR_{M,k}L_k^\top-X_k\cR_{M,k},
\end{equation}
where $X_k$ acts componentwise on matrices.

\begin{lemma}\label{lem:matrix}
Along the linear lifted evolution $\partial_t F=\sum_{k=1}^3 X_k^2 F$, any smooth symmetric matrix weight $M(v,w)$ satisfies
\begin{equation}\label{eq:matrix-identity}
 \frac{\dd}{\dd t}\iint_{\R^6}F \tilde s^\top M\tilde s\dd v\dd w
 =-2\sum_{k=1}^3\iint_{\R^6}FV_k^\top MV_k\dd v\dd w
   +\sum_{k=1}^3\iint_{\R^6}F\tilde s^\top\cK_{M,k}\tilde s\dd v\dd w.
\end{equation}
Evaluating the directional derivative \eqref{eq:matrix-identity} at the product density $F=f\otimes f$ with $M=S$ along the Landau equation \eqref{eq:landau} yields
\begin{equation}\label{eq:physical-factor}
 \frac{\dd}{\dd t}I_{-1}(f)=\frac{1}{16\pi}
 \left\{-2\sum_{k=1}^3\iint_{\R^6}FV_k^\top SV_k\dd v\dd w
       +\sum_{k=1}^3\iint_{\R^6}F\tilde s^\top\cK_{S,k}\tilde s\dd v\dd w\right\}.
\end{equation}
\end{lemma}

\begin{proof}
By linearity in the generator, it suffices to verify \eqref{eq:matrix-identity} for a single divergence-free vector field $X=X_k$, writing $x=X\log F$, $V=\nabla_{v,w}x$, and $L=D\Xi_k$. Along the flow $\partial_t F=X^2 F$, the logarithmic derivative satisfies
\[
 \partial_t\log F=Xx+x^2,\qquad V=X\tilde s+L^\top\tilde s.
\]
Differentiating the weighted Fisher integral in time and integrating by parts along the divergence-free field $X$ gives
\begin{align*}
 \frac{\dd}{\dd t}\iint_{\R^6}F\tilde s^\top M\tilde s\dd v\dd w
 &=-2\iint_{\R^6}FV^\top MV\dd v\dd w
   +2\iint_{\R^6}F\tilde s^\top\cR_M V\dd v\dd w\\
 &\quad +\iint_{\R^6}Fx\,\tilde s^\top\cR_M\tilde s\dd v\dd w.
\end{align*}
Using the identity $Fx=XF$ and integrating by parts along $X$ once more in the last term, we obtain
\[
 \iint_{\R^6}Fx\,\tilde s^\top\cR_M\tilde s\dd v\dd w
 =-2\iint_{\R^6}F(X\tilde s)^\top\cR_M\tilde s\dd v\dd w
   -\iint_{\R^6}F\tilde s^\top(X\cR_M)\tilde s\dd v\dd w.
\]
Substituting $X\tilde s=V-L^\top\tilde s$ into the first integral on the right-hand side and collecting terms according to \eqref{eq:RK} establishes \eqref{eq:matrix-identity}.

To deduce \eqref{eq:physical-factor}, we compute the first variation of the functional $\iint_{\R^6}F\tilde s^\top S\tilde s\dd v\dd w$ at $F=f\otimes f$. Setting
\[
 \Phi_\omega(v)=-2\diver_v\bigl(\omega(v)\nabla\log f(v)\bigr)
                -\omega(v)|\nabla\log f(v)|^2,
\]
the first variation with respect to $F$ evaluated at $f\otimes f$ is the sum $\Phi_\omega(v)+\Phi_\omega(w)$. Integrating the lifted collision operator $\sum_{k=1}^3 X_k^2(f\otimes f)$ with respect to $w$ (or $v$) recovers $8\pi Q(f,f)(v)$ (or $8\pi Q(f,f)(w)$). By symmetry of the two marginals, the directional derivative of the lifted functional along $\sum_{k=1}^3 X_k^2(f\otimes f)$ equals
\[
 16\pi\int_{\R^3}\Phi_\omega(v)Q(f,f)(v)\dd v
 =16\pi\,\frac{\dd}{\dd t}I_{-1}(f),
\]
which gives \eqref{eq:physical-factor}.
\end{proof}

We next introduce the pair dissipation $D_\Omega$ and the spherical gradient integral $R_\Omega$:
\begin{equation}\label{eq:D-R}
 D_\Omega=\sum_{k=1}^3\iint_{\R^6}F\Omega|V_k|^2\dd v\dd w,\qquad
 R_\Omega=\sum_{k=1}^3\iint_{\R^6}F\Omega\frac{x_k^2}{r^2}\dd v\dd w.
\end{equation}
When $M=I_6$, we factor $\Xi_k=\beta B_k$ with $\beta(r)=r^{-3/2}$ and $B_k(z)=(b_k(z),-b_k(z))$. Because $B_k$ is an infinitesimal rotation in $z$, its Jacobian $DB_k$ is skew-symmetric and orthogonal to $\nabla_{v,w}\beta$. Differentiating this orthogonality relation and using the skew-symmetry of $DB_k$ yields
\[
 B_k\cdot\nabla_{v,w}\beta=0,\qquad
 (DB_k)\nabla_{v,w}\beta=(D^2_{v,w}\beta)B_k,\qquad
 |\nabla_{v,w}\beta|^2=\frac{9}{2}r^{-5}.
\]
Substituting these relations into \eqref{eq:RK}, we find that the second derivative of $\beta$ cancels in $\cK_{I_6,k}$, leaving
\[
 \cR_{I_6,k}=B_k\otimes\nabla_{v,w}\beta+\nabla_{v,w}\beta\otimes B_k,\qquad
 \cK_{I_6,k}=2|\nabla_{v,w}\beta|^2B_k\otimes B_k
 =9r^{-2}\Xi_k\otimes\Xi_k.
\]
Contracting $\cK_{I_6,k}$ with the score $\tilde s$ produces the spherical curvature identity
\begin{equation}\label{eq:identity-I}
 \tilde s^\top\cK_{I_6,k}\tilde s=9\,\frac{x_k^2}{r^2}.
\end{equation}
Since $X_k\Omega=0$, the derivation rule in \eqref{eq:RK} immediately implies
\begin{equation}\label{eq:identity-Omega}
 \cK_{\Omega I_6,k}=\Omega\cK_{I_6,k}.
\end{equation}

Combining Ji's sharp spectral gap on the projective sphere \cite{JiSphere} with the comparison between $S$ and $\Omega I_6$ yields the baseline coercivity of the leading dissipation term in \eqref{eq:physical-factor}.

\begin{lemma}\label{lem:sphere}
For $F=f\otimes f$, the functionals in \eqref{eq:D-R} satisfy
\begin{equation}\label{eq:sphere}
 D_\Omega\ge 11R_\Omega.
\end{equation}
Furthermore,
\begin{equation}\label{eq:baseline}
 -2\sum_{k=1}^3\iint_{\R^6}FV_k^\top SV_k\dd v\dd w+9R_\Omega
 \le -\kappa_0 D_\Omega,\qquad
 \kappa_0=\sqrt{2}-\frac{9}{11}>0.
\end{equation}
\end{lemma}

\begin{proof}
On each collision sphere of fixed center of mass $c$ and radius $r$, the restricted product density $n\mapsto F(c+rn/2,c-rn/2)$ is an even function on $S^2$ and therefore descends to the real projective plane $\mathbb{RP}^2$. By Ji's Bakry--\'Emery $\Gamma_2$ inequality on $\mathbb{RP}^2$ \cite[Theorem 1.1]{JiSphere},
\[
 \int_{S^2}F\bigl(|D^2_{S^2}\log F|^2+|\nabla_{S^2}\log F|^2\bigr)\dd\sigma(n)
 \ge \frac{11}{2}\int_{S^2}F|\nabla_{S^2}\log F|^2\dd\sigma(n).
\]
In terms of the angular momentum operators $\mathcal L_k=(e_k\times n)\cdot\nabla_{S^2}$ on $S^2$, changing variables from $(v,w)$ to $(c,z)$ and discarding the center-of-mass and radial components of $V_k=\nabla_{v,w}x_k$ give
\[
 x_k=2r^{-3/2}\mathcal L_k\log F,\qquad
 |V_k|^2=\frac{1}{2}|\nabla_c x_k|^2+2|\nabla_z x_k|^2\ge 2r^{-2}|\nabla_{S^2}x_k|^2,
\]
together with the spherical Bochner identities
\[
 \sum_{k=1}^3|\mathcal L_k\log F|^2=|\nabla_{S^2}\log F|^2,\qquad
 \sum_{k=1}^3|\nabla_{S^2}\mathcal L_k\log F|^2=|D^2_{S^2}\log F|^2+|\nabla_{S^2}\log F|^2.
\]
Summing over $k\in\{1,2,3\}$ on each sphere yields the lower bound
\[
 \sum_{k=1}^3|V_k|^2\ge 8r^{-5}\bigl(|D^2_{S^2}\log F|^2+|\nabla_{S^2}\log F|^2\bigr),\qquad
 \sum_{k=1}^3\frac{x_k^2}{r^2}=4r^{-5}|\nabla_{S^2}\log F|^2.
\]
Since the pair weight $\Omega(v,w)$ depends only on $(c,r)$, multiplying by $F\Omega$ and integrating over $S^2$ and $(c,r)$ establishes \eqref{eq:sphere}.

Finally, for every $(v,w)\in\R^6$, the ratio of $\omega(v)$ to $\Omega(v,w)$ satisfies
\[
 \frac{\omega(v)}{\Omega(v,w)}
 =\left(\frac{2+|v|^2+|w|^2}{2(1+|v|^2)}\right)^{1/2}
 \ge \frac{1}{\sqrt{2}},
\]
and the same bound holds for $\omega(w)/\Omega(v,w)$ by symmetry. Thus $S\ge 2^{-1/2}\Omega I_6$ in the sense of symmetric matrices, which combines with \eqref{eq:sphere} to give \eqref{eq:baseline}.
\end{proof}

\section{The weight commutator}\label{sec:commutator}

Decomposing the weight matrix as $S=\Omega I_6+T$ with $T=S-\Omega I_6$, we see from \eqref{eq:identity-I} and \eqref{eq:identity-Omega} that the commutator sum in \eqref{eq:physical-factor} splits into $9R_\Omega+E_T$, where
\[
 E_T=\sum_{k=1}^3\iint_{\R^6}F\tilde s^\top\cK_{T,k}\tilde s\dd v\dd w.
\]
To estimate $E_T$, fix a smooth nonincreasing cutoff $\theta:[0,\infty)\to[0,1]$ equal to $1$ on $[0,1]$ and $0$ on $[2,\infty)$, and define the radial cutoff
\[
 \chi_\delta(v,w)=\theta\left(\frac{r}{\delta\rho}\right),\qquad 0<\delta<\frac{1}{8}.
\]
Because $\chi_\delta$ depends only on $r$ and $\rho=\la{c}$, we have $X_k\chi_\delta=0$ and hence $\cK_{\chi_\delta T,k}=\chi_\delta\cK_{T,k}$. This splits $E_T=E_1+E_2$ into the integral $E_1$ over nearby pairs supported on $\{r\le 2\delta\rho\}$ and the integral $E_2$ over separated pairs supported on $\{r\ge\delta\rho\}$.

\subsection{Commutator estimate for nearby pairs}

In the region $\{r\le 2\delta\rho\}$ with $0<\delta<1/8$, the particle velocities are comparable to the center of mass:
\[
 \la{v}\sim\la{w}\sim\rho,\qquad \Omega\sim\rho^{-1}.
\]
Since the perturbation matrix $T(v,w)=S(v,w)-\Omega(v,w)I_6$ vanishes identically on the collision diagonal $\{z=0\}$, Taylor's theorem around $z=0$ and the explicit derivatives of $\omega$ and $\Omega$ yield, uniformly for $0<\delta<1/8$,
\begin{equation}\label{eq:near-coefficients}
 \begin{aligned}
 |T|&\le C\frac{\Omega r}{\rho},&
 |DT|&\le C\frac{\Omega}{\rho},&
 |D^2 T|&\le C\frac{\Omega}{\rho^2},\\
 |\cR_{T,k}|&\le C r^{-3/2}\frac{\Omega r}{\rho},&
 |\diver_{v,w}(\chi_\delta\cR_{T,k})|&\le C r^{-3/2}\frac{\Omega}{\rho}.\span\omit
 \end{aligned}
\end{equation}
On the cutoff transition annulus $\{\delta\rho\le r\le 2\delta\rho\}$, we have
\[
 |\nabla_{v,w}\chi_\delta|\le C(\delta\rho)^{-1},\qquad
 |\cR_{T,k}|\le Cr^{-3/2}\frac{\Omega r}{\rho}\le 2C\delta r^{-3/2}\Omega,
\]
so the constant $C$ in the divergence estimate of \eqref{eq:near-coefficients} is independent of $\delta$.

If one expanded the second commutator $\cK_{T,k}$ directly at small relative velocities, differentiating $L_k=O(r^{-3/2})$ would produce a singular $r^{-2}\Omega/\rho$ kernel that cannot be controlled by the first-order Fisher information without higher moments. Instead, we undo one rotational derivative and integrate by parts once in Euclidean coordinates to transfer a derivative onto the logarithmic Hessian $D^2_{v,w}\log F$:
\begin{equation}\label{eq:near-identity}
 \begin{aligned}
 \iint_{\R^6}F\chi_\delta\tilde s^\top\cK_{T,k}\tilde s\dd v\dd w
 &=\iint_{\R^6}F\chi_\delta V_k^\top\cR_{T,k}\tilde s\dd v\dd w
   -\iint_{\R^6}Fx_k\,\diver_{v,w}(\chi_\delta\cR_{T,k})\cdot\tilde s\dd v\dd w\\
 &\quad-\iint_{\R^6}F\chi_\delta x_k\,\cR_{T,k}:D^2_{v,w}\log F\dd v\dd w.
 \end{aligned}
\end{equation}
To verify \eqref{eq:near-identity}, we first reverse the rotational integration by parts along $X_k$ from the proof of Lemma~\ref{lem:matrix}, which rewrites the left-hand side of \eqref{eq:near-identity} as
\[
 2\iint_{\R^6}F\chi_\delta\tilde s^\top\cR_{T,k}V_k\dd v\dd w
 +\iint_{\R^6}F\chi_\delta x_k\,\tilde s^\top\cR_{T,k}\tilde s\dd v\dd w.
\]
Writing one score factor in the second integral as $F\tilde s=\nabla_{v,w}F$, integrating by parts in $(v,w)$, and using $\nabla_{v,w}x_k=V_k$ then gives \eqref{eq:near-identity}.

Let $H_v=D^2\log f(v)$ denote the logarithmic Hessian of $f$, so that $D^2_{v,w}\log F=\diag(H_v,H_w)$. Using the coefficient bounds \eqref{eq:near-coefficients} and Cauchy--Schwarz with a small parameter $\varepsilon>0$, we bound the three integrands on the right-hand side of \eqref{eq:near-identity} by
\[
 \begin{aligned}
 |\chi_\delta V_k^\top\cR_{T,k}\tilde s|
 &\le\varepsilon\Omega|V_k|^2
   +C_\varepsilon\frac{\Omega}{\rho^2r}|\tilde s|^2,\\
 |x_k\,\diver_{v,w}(\chi_\delta\cR_{T,k})\cdot\tilde s|
 &\le\varepsilon\Omega\frac{x_k^2}{r^2}
   +C_\varepsilon\frac{\Omega}{\rho^2r}|\tilde s|^2,\\
 |\chi_\delta x_k\,\cR_{T,k}:D^2_{v,w}\log F|
 &\le\varepsilon\Omega\frac{x_k^2}{r^2}
   +C_\varepsilon\frac{\Omega r}{\rho^2}\bigl(|H_v|^2+|H_w|^2\bigr).
 \end{aligned}
\]
Multiplying by $F$, integrating over $\R^6$, and summing over $k\in\{1,2,3\}$ yields
\begin{equation}\label{eq:near-bound}
 E_1\le \varepsilon D_\Omega+2\varepsilon R_\Omega+C_\varepsilon P+C_\varepsilon Q,
\end{equation}
where the score remainder $P$ and Hessian remainder $Q$ are defined by
\begin{equation}\label{eq:PQ}
 \begin{aligned}
 P&=\iint_{\{r\le 2\delta\rho\}}f(v)f(w)\frac{\Omega}{\rho^2 r}\bigl(|s_v|^2+|s_w|^2\bigr)\dd v\dd w,\\
 Q&=\iint_{\{r\le 2\delta\rho\}}f(v)f(w)\frac{\Omega r}{\rho^2}\bigl(|H_v|^2+|H_w|^2\bigr)\dd v\dd w.
 \end{aligned}
\end{equation}

To bound the score integral $P$, we introduce the weighted density $\tilde f(w)=\omega(w)f(w)$. On $\R^3$, the translated Hardy inequality and Cauchy--Schwarz give
\begin{equation}\label{eq:Hardy}
 \int_{\R^3}\frac{\tilde f(w)}{|v-w|^2}\dd w\le I(\tilde f),\qquad
 \bigl\||\cdot|^{-1}*\tilde f\bigr\|_{L^\infty(\R^3)}
 \le \sqrt{\|\tilde f\|_{L^1(\R^3)}I(\tilde f)}.
\end{equation}
Since $|\nabla\log\omega(w)|\le\la{w}^{-1}\le 1$, expanding $\nabla\sqrt{\tilde f}$ gives $I(\tilde f)\le 2I_{-1}(f)+C$. On the domain $\{r\le 2\delta\rho\}$, we have $f(w)\le C\la{v}\tilde f(w)$ and $\Omega/\rho^2\le C\la{v}^{-3}$. By symmetry in $(v,w)$ and \eqref{eq:Hardy}, we obtain
\begin{equation}\label{eq:P}
 P\le C\int_{\R^3}f(v)\la{v}^{-2}|s_v|^2\left(\int_{\R^3}\frac{\tilde f(w)}{|v-w|}\dd w\right)\dd v
 \le C I_{-1}(f)\sqrt{1+I_{-1}(f)}.
\end{equation}

For the Hessian remainder $Q$, velocity comparability $\la{w}\sim\la{v}$ on $\{r\le 2\delta\rho\}$ and the conserved energy yield quadratic local mass decay at large velocities:
\[
 \int_{\{r\le 2\delta\rho\}}f(w)\dd w
 \le C\la{v}^{-2}\int_{\R^3}\la{w}^2 f(w)\dd w
 \le C\la{v}^{-2},
\]
and consequently, since $r\le 2\delta\rho\le C\delta\la{v}$,
\[
 \int_{\{r\le 2\delta\rho\}}rf(w)\dd w\le C\delta\la{v}^{-1}.
\]
Combining this linear gain in $\delta$ and decay in $\la{v}^{-1}$ with the kernel factor $\Omega/\rho^2\le C\la{v}^{-3}$, and using symmetry in $(v,w)$, we arrive at
\begin{equation}\label{eq:Q}
 Q\le C\delta\int_{\R^3}f(v)\la{v}^{-4}|H_v|^2\dd v
 =C\delta J_4(f).
\end{equation}

\subsection{Commutator estimate for separated pairs}

On the complementary region $\{r\ge\delta\rho\}$, the relative velocity $r$ is bounded away from zero and dominates both particle velocities:
\begin{equation}\label{eq:far-geometry}
 r\ge \delta,\qquad \max\bigl(\la{v},\la{w}\bigr)\le C_\delta r.
\end{equation}
Expanding $\cK_{T,k}=\cK_{S,k}-\Omega\cK_{I_6,k}$ directly from \eqref{eq:RK} and using the bounds $|L_k|\le Cr^{-3/2}$ and $|X_k L_k|\le Cr^{-3}$, we have
\begin{equation}\label{eq:far-K}
 \begin{aligned}
 |\cK_{T,k}|
 &\le C\Bigl\{r^{-3}\bigl(\omega(v)+\omega(w)+\Omega\bigr)
   +r^{-2}\bigl(|\nabla\omega(v)|+|\nabla\omega(w)|\bigr)\\
 &\qquad\quad +r^{-1}\bigl(|D^2\omega(v)|+|D^2\omega(w)|\bigr)\Bigr\}.
 \end{aligned}
\end{equation}
Since $\omega(v)=\la{v}^{-1}$ and its first two derivatives are globally bounded on $\R^3$, the geometric comparability \eqref{eq:far-geometry} implies
\[
 |\cK_{T,k}(v,w)|\le C_\delta r^{-1}\le C_\delta\min\bigl(\omega(v),\omega(w)\bigr)\qquad\text{on }\{r\ge\delta\rho\}.
\]
Integrating against $F|\tilde s|^2=f(v)f(w)(|s_v|^2+|s_w|^2)$ and using $\int_{\R^3}f\dd v=1$ therefore gives
\begin{equation}\label{eq:far-bound}
 |E_2|\le C_\delta I_{-1}(f).
\end{equation}

Combining the evolution identity \eqref{eq:physical-factor}, the spherical coercivity \eqref{eq:baseline}, the near-diagonal bounds \eqref{eq:near-bound}, \eqref{eq:P}, and \eqref{eq:Q}, and the separated-pair bound \eqref{eq:far-bound} with $\varepsilon=\kappa_0/4$, we obtain
\begin{equation}\label{eq:before-coercivity}
 \frac{\dd}{\dd t}I_{-1}(f)
 \le -c_0 D_\Omega+C\delta J_4(f)+C_\delta\bigl(1+I_{-1}(f)\bigr)^{3/2},
\end{equation}
where $c_0=\kappa_0/(32\pi)>0$.

\section{Coercivity, cutoff justification, and proof of Theorem~\ref{thm:critical-fisher}}\label{sec:coercivity}

\subsection{Coercivity of the pair dissipation and Riccati closure}

To turn \eqref{eq:before-coercivity} into a closed differential inequality for $I_{-1}(f)$, we next prove that the pair dissipation $D_\Omega$ controls the weighted logarithmic Hessian functional $J_4(f)$.

\begin{lemma}\label{lem:coercivity}
Under the normalization \eqref{eq:normalization}, there exist constants $c_H,C_H>0$ depending only on $H$ such that
\begin{equation}\label{eq:coercivity}
 D_\Omega\ge c_H J_4(f)-C_H.
\end{equation}
\end{lemma}

\begin{proof}
To regularize the $r^{-3/2}$ singularity at small relative velocities, we insert the radial multiplier $0\le\eta(r)=\min(1,r^3)\le 1$ into $D_\Omega$. Writing $\tilde b_k(z)=r^{-3/2}b_k(z)$, the $v$-component of $V_k=\nabla_{v,w}x_k$ expands as
\[
 \nabla_v x_k=H_v\tilde b_k+(D\tilde b_k)^\top(s_v-s_w),\qquad
 \sum_{k=1}^3\bigl|(D\tilde b_k)^\top q\bigr|^2\le C r^{-3}|q|^2\quad (q\in\R^3).
\]
Using the elementary inequality $|\nabla_v x_k|^2\ge\frac{1}{2}|H_v\tilde b_k|^2-|(D\tilde b_k)^\top(s_v-s_w)|^2$, together with the identical bound for $\nabla_w x_k$ and the identity $\sum_{k=1}^3 b_k(z)\otimes b_k(z)=r^2(I_3-\tilde z\otimes\tilde z)$, we obtain
\begin{equation}\label{eq:coercive-expansion}
 D_\Omega\ge c\int_{\R^3}f(v)\,\tr\bigl(H_v\mathcal B(v)H_v^\top\bigr)\dd v
   -C\iint_{\R^6}f(v)f(w)\Omega\min(1,r^{-3})\bigl(|s_v|^2+|s_w|^2\bigr)\dd v\dd w,
\end{equation}
where the $3\times 3$ matrix field $\mathcal B(v)$ is given by
\[
 \mathcal B(v)=\int_{\R^3}f(w)\Omega(v,w)\frac{\eta(r)}{r}\bigl(I_3-\tilde z\otimes\tilde z\bigr)\dd w.
\]

We claim that $\mathcal B(v)\ge c_H\la{v}^{-4}I_3$ for all $v\in\R^3$. By the energy normalization $\int_{\R^3}|w|^2f(w)\dd w=3$, the negative part of $f\log f$ on $\{f\le 1\}$ is controlled by splitting at $f=e^{-1-|w|^2}$ and using $x|\log x|\le C\sqrt{x}$ on $(0,1]$, which gives $\int_{\R^3}f\log_+ f\dd w\le C_H$. Consequently, for every measurable set $E\subset\R^3$ and every threshold $M>1$,
\[
 \int_E f(w)\dd w\le M|E|+\frac{C_H}{\log M}.
\]
By Chebyshev's inequality, the ball $B_{\sqrt{12}}$ carries mass at least $3/4$. For any unit vector $\xi\in S^2$ and any $v\in\R^3$, the cylinder of radius $d>0$ around the line $v+\R\xi$ intersects $B_{\sqrt{12}}$ in a set of Lebesgue measure at most $Cd^2$. We choose $M=M(H)>1$ and $d=d_H\in(0,1)$ so that
\[
 \frac{C_H}{\log M}\le\frac{1}{8},\qquad CMd_H^2\le\frac{1}{8}.
\]
Then the cylinder $\{w\in B_{\sqrt{12}}:|(v-w)\times\xi|<d_H\}$ carries mass at most $1/4$. On the complementary subset of $B_{\sqrt{12}}$, which has mass at least $1/2$, we have
\[
 |(v-w)\times\xi|\ge d_H,\qquad d_H\le r\le C\la{v},\qquad \Omega(v,w)\ge c\la{v}^{-1},
\]
and therefore
\[
 \Omega(v,w)\frac{\eta(r)}{r}|\tilde z\times\xi|^2
 =\Omega(v,w)\frac{\min(1,r^3)}{r^3}|(v-w)\times\xi|^2
 \ge c_H\la{v}^{-4}.
\]
Integrating against $f(w)\dd w$ over this subset proves $\xi^\top\mathcal B(v)\xi\ge c_H\la{v}^{-4}$.

To bound the first-order score error in \eqref{eq:coercive-expansion}, we verify that for all $v\in\R^3$,
\begin{equation}\label{eq:error-convolution}
 \int_{\R^3}f(w)\Omega(v,w)\min(1,r^{-3})\dd w\le C\la{v}^{-3}.
\end{equation}
This is immediate for $|v|\le 2$. For $|v|>2$, we split the integral at $r=\la{v}/2$: on $\{r\le\la{v}/2\}$ we have $\la{w}\ge\la{v}/2$, so $\Omega\le C\la{v}^{-1}$ and Chebyshev's inequality with the conserved energy gives
\[
 \int_{\{r\le\la{v}/2\}}f(w)\dd w\le C\la{v}^{-2},
\]
whereas on $\{r>\la{v}/2\}$ we simply have $\min(1,r^{-3})\le C\la{v}^{-3}$. Substituting \eqref{eq:error-convolution} into \eqref{eq:coercive-expansion} yields
\[
 D_\Omega\ge 2c_H J_4(f)-CI_{-3}(f).
\]
Finally, integrating by parts in $I_{-3}(f)$ and using $|\Delta(\la{v}^{-3})|\le C\la{v}^{-5}$ give
\[
 I_{-3}(f)=-\int_{\R^3}f\la{v}^{-3}\Delta\log f\dd v+\int_{\R^3}f\Delta\bigl(\la{v}^{-3}\bigr)\dd v
 \le \sqrt{3J_4(f)}+C\le c_H J_4(f)+C_H,
\]
which absorbs into $2c_H J_4(f)$ and completes the proof of \eqref{eq:coercivity}.
\end{proof}

Choosing $\delta=\delta(H)>0$ small enough in \eqref{eq:before-coercivity} that $C\delta\le\frac{1}{2}c_0c_H$ and applying Lemma~\ref{lem:coercivity} establishes the first inequality in \eqref{eq:main-diff}. To close the Riccati inequality using only the conserved second moment $\int_{\R^3}f\la{v}^2\dd v=4$, we integrate $I_{-1}(f)$ by parts with $\omega(v)=\la{v}^{-1}$ and use the identity $\Delta\omega=-3\la{v}^{-5}<0$:
\begin{equation}\label{eq:exact-J4}
 I_{-1}(f)+3\int_{\R^3}f\la{v}^{-5}\dd v
 =-\int_{\R^3}\omega f\Delta\log f\dd v
 \le \left(3J_4(f)\int_{\R^3}f\la{v}^2\dd v\right)^{1/2}
 =\sqrt{12J_4(f)}.
\end{equation}
Dropping the nonnegative second term on the left-hand side gives $J_4(f)\ge\frac{1}{12}I_{-1}(f)^2$. Absorbing the lower-order power $C_H(1+I_{-1}(f))^{3/2}$ into $-\frac{1}{24}c_H I_{-1}(f)^2$ via Young's inequality proves the second inequality in \eqref{eq:main-diff}. Once the local absolute continuity of $t\mapsto I_{-1}(f(t))$ on $(0,T]$ is justified in Section~\ref{app:justification}, comparing $y(t)=I_{-1}(f(t))$ on $(s,T]$ with the explicit supersolution
\[
 Y_s(t)=\left(\frac{\tilde C_H}{\tilde c_H}\right)^{1/2}+\frac{1}{\tilde c_H(t-s)},\qquad t>s>0,
\]
yields \eqref{eq:main-bound} after letting $s\downarrow 0$.

\subsection{Approximation and cutoff justification}\label{app:justification}

We now justify the preceding formal calculations for any given smooth, rapidly decaying solution $f$ on $[0,T]\times\R^3$. Qualitative bounds depending on the Schwartz seminorms of $f$ are used solely to ensure integrability and pass to the limit in spatial cutoffs; they do not enter the quantitative constants in Theorem~\ref{thm:critical-fisher}.

First, for the initial datum $f_0=f(0,\cdot)$, Taylor's theorem along each coordinate axis together with $f_0\ge 0$ gives the classical bound
\[
 \frac{|\nabla f_0(v)|^2}{f_0(v)}
 \le C f_0(v)+C\sup_{u\in B_1(v)}|D^2f_0(u)|,
\]
which is integrable on $\R^3$ because $f_0$ is rapidly decaying. Thus $I(f_0)<\infty$ and $\sqrt{f_0}\in H^1(\R^3)$. Setting
\[
 D_0(f)=\sum_{k=1}^3\iint_{\R^6}F|V_k|^2\dd v\dd w,
 \qquad c_1=\frac{2-9/11}{16\pi}=\frac{13}{176\pi},
\]
we first establish the unweighted dissipation bound starting from $t=0$:
\begin{equation}\label{eq:integrated-ordinary-fisher}
 I(f(t))+c_1\int_0^tD_0(f(s))\dd s\le I(f_0),
 \qquad 0<t\le T.
\end{equation}

To prove \eqref{eq:integrated-ordinary-fisher}, fix a radial cutoff $0\le\chi\le1$ in $C_c^\infty(\R^3)$ equal to $1$ on $B_1$, and define the regularized initial data
\[
 g_{0,j}(v)=\chi(v/j)^2f_0(v)+\varepsilon_j e^{-\beta|v|^2},
 \qquad
 m_j=\int_{\R^3} g_{0,j}(v)\dd v,\qquad f_{0,j}=m_j^{-1}g_{0,j},
\]
with $\beta>0$ and $\varepsilon_j\downarrow 0$. Then $m_j\to 1$, each $f_{0,j}$ is a smooth probability density bounded above and below by positive multiples of $e^{-\beta|v|^2}$, and $f_{0,j}\to f_0$ in every polynomially weighted $L^1$ space. Since $\chi(\cdot/j)\sqrt{f_0}\to\sqrt{f_0}$ in $H^1(\R^3)$, the convexity and homogeneity of the Fisher information combined with lower semicontinuity give
\[
 \limsup_{j\to\infty}I(f_{0,j})
 \le\lim_{j\to\infty}\frac{1}{m_j}
 \bigl[I\bigl(\chi(\cdot/j)^2f_0\bigr)
   +\varepsilon_jI(e^{-\beta|v|^2})\bigr]=I(f_0),
\]
so that $I(f_{0,j})\to I(f_0)$ as $j\to\infty$.

Let $f_j$ denote the global smooth solution with initial datum $f_{0,j}$ constructed by Guillen and Silvestre \cite[Theorem 1.2 and Appendix B]{GS}. Because $f_j$ enjoys two-sided Gaussian bounds and polynomial derivative bounds on $[0,T]\times\R^3$, the unweighted Fisher identity and Lemma~\ref{lem:sphere} with weight $\Omega\equiv 1$ apply directly to $f_j$, giving
\[
 I(f_j(t))+c_1\int_0^tD_0(f_j(s))\dd s\le I(f_{0,j}).
\]
By Sobolev embedding for $\sqrt{f_j}$ and polynomial moment propagation \cite{Desvillettes}, the sequence $\{f_j\}$ has uniformly bounded moments of all orders on $[0,T]$. Interpolation then gives, for any fixed $\kappa\ge 5$,
\[
 \int_{\R^3}\la v^\kappa f_j(t,v)^2\dd v
 \le\left(\int_{\R^3} f_j(t,v)^3\dd v\right)^{1/2}
     \left(\int_{\R^3}\la v^{2\kappa}f_j(t,v)\dd v\right)^{1/2}
 \le C_{f_0,T,\kappa},
\]
while Hardy's inequality and classical entropy ellipticity provide the uniform coefficient bounds
\[
 \|a[f_j]\|_{L^\infty}+\|\nabla a[f_j]\|_{L^\infty}\le C_{f_0,T},\qquad
 c\la{v}^{-5}\le\det A[f_j](t,v),\qquad
 0\le A[f_j](t,v)\le\Lambda I_3.
\]
Applying Silvestre's conditional parabolic regularization theorem \cite[Theorems 1.1 and 2.1]{Silvestre} (with dimension $3$, $p=2$, $\beta=0$, and $\kappa\ge 5$) yields the uniform $L^\infty$ bound
\begin{equation}\label{eq:approximation-uniform-bound}
 \|f_j(t)\|_{L^\infty(\R^3)}
 \le C_{f_0,T}\bigl(1+t^{-3/4}\bigr),\qquad 0<t\le T.
\end{equation}

On compact subsets of $(0,T]\times\R^3$, interior De Giorgi--Nash--Moser and parabolic Schauder estimates provide uniform $C^2_{\mathrm{loc}}$ compactness for $\{f_j\}$. Moreover, for any test function $\varphi\in C_c^2(\R^3)$, the weak formulation
\[
 \frac{\dd}{\dd t}\int_{\R^3}\varphi f_j\dd v
 =\int_{\R^3} f_j\bigl(A[f_j]:D^2\varphi+2\nabla a[f_j]\cdot\nabla\varphi\bigr)\dd v
\]
shows that the time derivative of $\int_{\R^3}\varphi f_j\dd v$ is uniformly bounded on $[0,T]$. Passing to a subsequence, $f_j$ converges in $C^2_{\mathrm{loc}}((0,T]\times\R^3)$ to a weak solution in $L^1(0,T;L^\infty(\R^3))$ with uniformly bounded second moment and initial trace $f_0$. Since $t^{-3/4}$ is integrable at $t=0$, Fournier's uniqueness theorem for bounded solutions \cite{Fournier} identifies the limit with our given solution $f$. By the strong parabolic maximum principle, $f$ is strictly positive on $(0,T]\times\R^3$, so the nonnegative integrands of $I(f_j)$ and $D_0(f_j)$ converge almost everywhere on $\{v\ne w\}$. Fatou's lemma together with $I(f_{0,j})\to I(f_0)$ therefore proves \eqref{eq:integrated-ordinary-fisher}.

From \eqref{eq:integrated-ordinary-fisher} we have $I(f)\in L^\infty(0,T)$ and $D_\Omega(f)\le D_0(f)\in L^1(0,T)$. Applying \eqref{eq:coercive-expansion} and \eqref{eq:error-convolution} on compactly truncated domains and passing to the limit via Fatou's lemma gives
\[
 c_HJ_4(f(t))\le D_\Omega(f(t))+C I(f(t))
 \le D_\Omega(f(t))+C I(f_0)
\]
for almost every $t\in(0,T)$, so that $J_4(f)\in L^1(0,T)$. This immediately justifies the single-particle integrations by parts in Lemma~\ref{lem:coercivity} and \eqref{eq:exact-J4}: multiplying by a smooth radial cutoff $\vartheta_R(v)=\theta(|v|/R)$ gives, for $q\in\{1,3\}$,
\[
 \int_{\R^3}\vartheta_R\la v^{-q}f|\nabla\log f|^2\dd v
 =-\int_{\R^3}\vartheta_R\la v^{-q}f\Delta\log f\dd v
   +\int_{\R^3} f\Delta\bigl(\vartheta_R\la v^{-q}\bigr)\dd v,
\]
and the cutoff error vanishes as $R\to\infty$ because
\[
 \la{v}^{-q}f|\Delta\log f|
 \le\frac{1}{2}f\la{v}^{-4}|D^2\log f|^2+\frac{3}{2}f\la{v}^{4-2q}\in L^1(\R^3)
\]
at almost every time. Similarly, \eqref{eq:sphere}, \eqref{eq:P}, and \eqref{eq:Q} imply that $D_\Omega,R_\Omega,P,Q\in L^1(0,T)$, while the translated Hardy inequality \eqref{eq:Hardy} and the near-diagonal bound $|\cK_{T,k}|\le C\Omega/(\rho r^2)$ ensure the space-time integrability
\[
 \iint_{\R^6}\frac{F}{r^2}\bigl(|s_v|^2+|s_w|^2\bigr)\dd v\dd w
 \le 2I(f)^2\le 2I(f_0)^2.
\]

Finally, to justify the lifted weighted evolution identity \eqref{eq:physical-factor} and the integration by parts \eqref{eq:near-identity} on any compact interval $[t_1,t_2]\subset(0,T]$, we introduce the truncated weights and cutoffs
\[
 \begin{aligned}
 \omega_N(v)&=\omega(v)\theta\left(\frac{|v|}{N}\right),&
 S_N(v,w)&=\diag\bigl(\omega_N(v)I_3,\omega_N(w)I_3\bigr),\\
 \zeta_L(c)&=\theta\left(\frac{|c|}{L}\right),&
 \eta_\varepsilon(r)&=1-\theta\left(\frac{r}{\varepsilon}\right).
 \end{aligned}
\]
Because $X_k\zeta_L=X_k\eta_\varepsilon=0$, applying the directional derivative identity of Lemma~\ref{lem:matrix} to the compactly supported smooth weight $\eta_\varepsilon\zeta_L S_N$ on $[t_1,t_2]\times\R^6$ reduces to a finite-dimensional integration by parts on a compact set where $f$ is smooth and strictly positive. On this compact cylinder, as $r\downarrow 0$,
\[
 \sum_{k=1}^3 X_k^2 F=O(r^{-1}),\qquad x_k=O(r^{1/2}),\qquad V_k=O(r^{-1/2}),
\]
so the boundary term involving $\nabla\eta_\varepsilon$ on $\{\varepsilon\le r\le 2\varepsilon\}$ has integral $O(\varepsilon)$ and vanishes as $\varepsilon\downarrow 0$. Letting $L\to\infty$ at fixed $N$ via the rapid decay of $f$, and then letting $N\to\infty$ by dominated convergence using the space-time integrable majorants $D_0$, $R_\Omega$, and $F(|s_v|^2+|s_w|^2)/r^2$, we conclude that $t\mapsto I_{-1}(f(t))$ is locally absolutely continuous on $(0,T]$ and satisfies \eqref{eq:physical-factor}. Inserting $\zeta_L(c)\eta_\varepsilon(r)$ into \eqref{eq:near-identity} and passing to the limit $\varepsilon\downarrow 0$ and $L\to\infty$ in the same manner completes the justification of \eqref{eq:near-bound} and finishes the proof of Theorem~\ref{thm:critical-fisher}.

\section{Global integrability, coefficient bounds, and local Moser iteration}\label{sec:local}

In this section we prove Corollary~\ref{cor:smoothing}. Setting $\tilde f=\omega f$ with $\omega(v)=\la{v}^{-1}$, the three-dimensional Sobolev inequality applied to $\sqrt{\tilde f}$ and the bound $I(\tilde f)\le 2I_{-1}(f)+C$ give
\begin{equation}\label{eq:weighted-L3}
 \|\omega f(t)\|_{L^3(\R^3)}\le C\bigl(1+I_{-1}(f(t))\bigr).
\end{equation}
Interpolating \eqref{eq:weighted-L3} with the conserved second moment via H\"older's inequality yields
\[
 \int_{\R^3}f^{9/5}\dd v
 \le \left(\int_{\R^3}\la{v}^{-3}f^3\dd v\right)^{2/5}\left(\int_{\R^3}\la{v}^2f\dd v\right)^{3/5},
\]
so that $\|f(t)\|_{L^{9/5}(\R^3)}\le C(1+I_{-1}(f(t)))^{2/3}$. Splitting the Newton potential $a[f]$ at radius $r_0>0$ and applying the Hardy--Littlewood--Sobolev inequality to $\nabla a[f]$ give
\[
 a[f](v)\le C\Bigl(\|f\|_{L^{9/5}(\R^3)}r_0^{1/3}+r_0^{-1}\Bigr),\qquad
 \|\nabla a[f]\|_{L^{9/2}(\R^3)}\le C\|f\|_{L^{9/5}(\R^3)}.
\]
Optimizing over $r_0>0$ yields $\|a[f]\|_{L^\infty(\R^3)}\le C\|f\|_{L^{9/5}(\R^3)}^{3/4}$, which combines with Theorem~\ref{thm:critical-fisher} to establish \eqref{eq:cor-bounds}.

\subsection{Local Moser iteration on bounded velocity balls}

Fix a radius $R\ge 1$ and $0<\tau<T$. On the cylinder $[\tau/2,T]\times B_{2R}$, the coefficient bounds \eqref{eq:cor-bounds} and classical entropy ellipticity imply
\[
 \lambda_R I_3\le A[f]\le \Lambda I_3,\qquad
 \lambda_R=c_H\la{2R}^{-3},\qquad
 \|b\|_{L^\infty([\tau/2,T];L^{9/2}(\R^3))}\le B,\qquad b=\nabla a[f].
\]
For any exponent $p\ge 9/5$, set $z=\vartheta\eta f^{p/2}$, where $\vartheta(t)$ and $\eta(v)$ are smooth nonnegative cutoffs supported in $[\tau/2,T]\times B_{2R}$. Testing the divergence-form equation $\partial_t f=\diver(A[f]\nabla f-bf)$ with $p\vartheta^2\eta^2 f^{p-1}$ gives
\begin{equation}\label{eq:moser-energy-pre}
 \frac{\dd}{\dd t}\int_{\R^3}z^2\dd v+\lambda_R\int_{\R^3}|\nabla z|^2\dd v
 \le C\int_{\R^3}\bigl(|\vartheta\vartheta'|\eta^2+\vartheta^2|\nabla\eta|^2\bigr)f^p\dd v
   +Cp^2\lambda_R^{-1}\int_{\R^3}|b|^2z^2\dd v.
\end{equation}
By H\"older's inequality, Gagliardo--Nirenberg interpolation in $\R^3$, and Young's inequality, the drift integral satisfies
\begin{equation}\label{eq:moser-drift}
 \int_{\R^3}|b|^2z^2\dd v
 \le B^2\|z\|_{L^{18/5}(\R^3)}^2
 \le CB^2\|\nabla z\|_{L^2(\R^3)}^{4/3}\|z\|_{L^2(\R^3)}^{2/3}
 \le \varepsilon\|\nabla z\|_{L^2(\R^3)}^2+C\varepsilon^{-2}B^6\|z\|_{L^2(\R^3)}^2.
\end{equation}
Choosing $\varepsilon=\frac{1}{2}c\lambda_R^2 p^{-2}$ absorbs the gradient term into the left-hand side of \eqref{eq:moser-energy-pre}. For each integer $j\ge 0$, define the nested space-time cylinders and iteration exponents
\[
 a_j=\tau(1-2^{-j-1}),\qquad
 r_j=R(1+2^{-j}),\qquad
 Q_j=[a_j,T]\times B_{r_j},\qquad
 p_j=\frac{9}{5}\left(\frac{5}{3}\right)^j.
\]
Choosing smooth cutoffs $\vartheta_j(t)$ and $\eta_j(v)$ between $Q_j$ and $Q_{j+1}$ with
\[
 |\vartheta_j'|\le \frac{C2^j}{\tau},\qquad
 |\nabla\eta_j|\le \frac{C2^j}{R},\qquad
 z_j=\vartheta_j\eta_j f^{p_j/2},
\]
and integrating \eqref{eq:moser-energy-pre} in time yields
\begin{equation}\label{eq:moser-energy}
 \sup_{a_j\le t\le T}\int_{\R^3}z_j(t,v)^2\dd v
 +\iint_{Q_j}|\nabla z_j|^2\dd v\dd t
 \le C_{H,\tau,R}p_j^6 4^j\bigl(1+\tau^{-1}+R^{-2}\bigr)\iint_{Q_j}f^{p_j}\dd v\dd t.
\end{equation}
Applying the three-dimensional parabolic Sobolev embedding on each cylinder gives the recurrence
\[
 \|f\|_{L^{5p_j/3}(Q_{j+1})}
 \le \Bigl[C_{H,\tau,R}p_j^6 4^j\bigl(1+\tau^{-1}+R^{-2}\bigr)\Bigr]^{1/p_j}\|f\|_{L^{p_j}(Q_j)}.
\]
Iterating over $j\ge 0$ from $p_0=9/5$ and passing to the limit $j\to\infty$ proves \eqref{eq:local-infty}.

\section{A homogeneous truncated mixed Nash inequality}\label{sec:truncated-mixed-nash}

The standard three-dimensional Nash inequality bounds the $L^{5/3}(\R^3)$ integral of a nonnegative density in terms of its mass and its full three-dimensional Fisher information. In the exterior velocity region of the Landau--Coulomb equation, however, the background diffusion matrix $A[f]$ degenerates like $\la{v}^{-3}$ in the radial direction and provides $\la{v}^{-1}$ ellipticity only in the two transverse directions. In this section we prove a nonlinear functional inequality showing that on any upper level set of the density, the Fisher information in only two transverse coordinate directions, when combined with the high-density part of the Coulomb entropy dissipation, suffices to control the $L^{5/3}$ norm of the truncated density without using any longitudinal derivative of $f$.

Write $v=(v_1,v_2,v_3)\in\R^3$ and $\nabla_\perp=(\partial_{v_2},\partial_{v_3})$. Let $f\ge 0$ be smooth and rapidly decaying on $\R^3$, and fix a density threshold $K>0$. We define the upper level set truncations and transverse Fisher energy by
\[
 \begin{aligned}
 H&=\{f>K\},&F&=f\mathbf{1}_H,&g&=(f-K)_+,&u&=(\sqrt{f}-\sqrt{K})_+,&h&=u^2,\\
 \mu&=\int_{\R^3}F\dd v,&X&=\int_H\frac{|\nabla_\perp f|^2}{f}\dd v=4\int_{\R^3}|\nabla_\perp u|^2\dd v,\span\omit\span\omit\span\omit
 \end{aligned}
\]
together with the Coulomb entropy dissipation restricted to $H\times H$:
\[
 \cD=\frac{1}{2}\iint_{H\times H}f(v)f(w)\bigl(\nabla\log f(v)-\nabla\log f(w)\bigr)^\top K(v-w)\bigl(\nabla\log f(v)-\nabla\log f(w)\bigr)\dd v\dd w.
\]

\begin{proposition}\label{prop:mixed-nash}
There exists a universal constant $C>0$, independent of $K$ and $f$, such that
\begin{equation}\label{eq:truncated-mixed-nash}
 \int_{\R^3}h^{5/3}\dd v
 \le C\left(\mu^{2/3}X+\mu^{1/2}X^{1/2}\cD^{1/3}\right).
\end{equation}
\end{proposition}

\begin{proof}
We may assume without loss of generality that $\mu>0$ and $X+\cD<\infty$; note that $X>0$ whenever $\mu>0$ and $X<\infty$, since otherwise $u\in L^2(\R^3)$ would be constant in $(v_2,v_3)$ and hence vanish identically.

\medskip
\noindent\emph{Step 1.}
Suppose $F_0\in L^1(\R^3)$ and $g_0\in W^{1,1}(\R^3)$ are nonnegative functions supported in a set of diameter at most $L>0$ and satisfying
\[
 0\le g_0\le F_0,\qquad F_0\sigma=\nabla g_0\quad\text{weakly on }\R^3,
\]
for a measurable vector field $\sigma=(\sigma_1,\sigma_2,\sigma_3)$ with $F_0\sigma\in L^1(\R^3)$ and $\sigma=0$ on $\{F_0=0\}$. Put
\[
 \mu_0=\int_{\R^3}F_0\dd v,\qquad b_0=\int_{\R^3}g_0\dd v>0,\qquad
 I_i=\int_{\R^3}F_0\sigma_i^2\dd v,\qquad X_0=I_2+I_3,
\]
and
\[
 \cD_0=\frac{1}{2}\iint_{\R^3\times\R^3}F_0(v)F_0(w)\bigl(\sigma(v)-\sigma(w)\bigr)^\top K(v-w)\bigl(\sigma(v)-\sigma(w)\bigr)\dd v\dd w.
\]
Define the centers of mass $\bar v$ and $\bar v_g$ and the covariance matrix $P_{ij}$ of $F_0$ by
\[
 \bar v=\frac{1}{\mu_0}\int_{\R^3}vF_0(v)\dd v,\qquad
 \bar v_g=\frac{1}{b_0}\int_{\R^3}vg_0(v)\dd v,\qquad
 P_{ij}=\int_{\R^3}F_0(v)(v_i-\bar v_i)(v_j-\bar v_j)\dd v.
\]
To reconstruct the longitudinal score component $\sigma_1(v)$ from the collision dissipation $\cD_0$ and the transverse components $\sigma_j$ ($j\in\{2,3\}$), we test the $(1,j)$-component of $(v-w)\times(\sigma(v)-\sigma(w))$ against $(w_j-\bar v_j)F_0(w)$ and integrate by parts in $w$ using the weak identity $F_0\sigma=\nabla g_0$. This yields the representation
\begin{equation}\label{eq:generalized-score-moment}
 \begin{split}
 &\int_{\R^3}(w_j-\bar v_j)F_0(w)\bigl[(v_1-w_1)\bigl(\sigma_j(v)-\sigma_j(w)\bigr)-(v_j-w_j)\bigl(\sigma_1(v)-\sigma_1(w)\bigr)\bigr]\dd w\\
 &\qquad =P_{jj}\sigma_1(v)-P_{1j}\sigma_j(v)+b_0(v_1-\bar v_{g,1}).
 \end{split}
\end{equation}
Integrating $F_0\sigma_j=\partial_{v_j}g_0$ by parts against $v_j-\bar v_j$ and using the support diameter bound $|v-\bar v|\le L$ give the variance bounds
\[
 P_{jj}I_j\ge b_0^2,\qquad P_{1j}^2\le P_{11}P_{jj},\qquad
 P_{11}\le L^2\mu_0,\qquad \int_{\R^3}F_0|v_1-\bar v_{g,1}|^2\dd v\le L^2\mu_0.
\]
Let $J_j(v)$ denote the left-hand side of \eqref{eq:generalized-score-moment}. Using the algebraic identity
\[
 |(v-w)\times(\sigma(v)-\sigma(w))|^2
 =8\pi|v-w|^3(\sigma(v)-\sigma(w))^\top K(v-w)(\sigma(v)-\sigma(w))
\]
together with $|v-w|\le L$ on the support of $F_0$, Cauchy--Schwarz in $w$ followed by integration against $F_0(v)\dd v$ gives
\[
 \int_{\R^3}F_0(v)|J_j(v)|^2\dd v\le CL^3P_{jj}\cD_0.
\]
Dividing \eqref{eq:generalized-score-moment} by $P_{jj}\ge b_0^2/I_j>0$, squaring, and integrating against $F_0(v)\dd v$, we obtain
\[
 I_1\le C\left(\frac{L^3\cD_0}{P_{jj}}
       +\frac{P_{1j}^2 I_j}{P_{jj}^2}
       +\frac{b_0^2 L^2\mu_0}{P_{jj}^2}\right)
 \le C\left(\frac{L^3 I_j\cD_0}{b_0^2}
       +\frac{L^2\mu_0 I_j^2}{b_0^2}\right).
\]
Since $I_j\le X_0$, this provides the cell longitudinal bound
\begin{equation}\label{eq:generalized-score-longitudinal}
 I_1\le C\left(\frac{L^3 X_0 \cD_0}{b_0^2}+\frac{L^2\mu_0 X_0^2}{b_0^2}\right).
\end{equation}

\medskip
\noindent\emph{Step 2.}
Fix a smooth quadratic partition of unity $\{\chi_k\}$ on $\R^3$ at scale $L>0$ with uniformly bounded overlap $N$, support diameter at most $CL$, and
\[
 \sum_k\chi_k^2=1,\qquad \sum_k|\nabla\chi_k|^2\le CL^{-2}.
\]
Define the localized cell functions and masses by
\[
 F_k=\chi_k^2 F,\quad g_k=\chi_k^2 g,\quad h_k=\chi_k^2 h,\quad
 \mu_k=\int_{\R^3}F_k\dd v,\quad b_k=\int_{\R^3}g_k\dd v,\quad n_k=\int_{\R^3}h_k\dd v,
\]
so that $0\le n_k\le b_k\le\mu_k$. On $H\cap\{\chi_k>0\}$, define the localized generalized score by
\begin{equation}\label{eq:truncated-generalized-score}
 \sigma_k=\frac{\nabla g_k}{F_k}=\nabla\log f+2\frac{g}{f}\nabla\log\chi_k,
\end{equation}
and set $\sigma_k=0$ where $F_k=0$. Because $\nabla g=\mathbf{1}_H\nabla f$ weakly on $\R^3$, the weak identity $F_k\sigma_k=\nabla g_k$ holds on $\R^3$ without any boundary measure on $\{f=K\}$. Let $I_{i,k},X_k,\cD_k$ denote the corresponding cell functionals from Step~1. On any cell with $b_k>0$, the gradient of $\sqrt{h_k}=\chi_k u$ differs from $\frac{1}{2}\sqrt{F_k}\,\sigma_k$ by the cutoff error
\begin{equation}\label{eq:truncated-gradient-error}
 \nabla(\chi_k u)-\frac{1}{2}\sqrt{F_k}\,\sigma_k=-\sqrt{\frac{K}{f}}\,u\nabla\chi_k.
\end{equation}
Setting
\[
 e_k=\int_{\R^3}h|\nabla\chi_k|^2\dd v,\qquad
 Z_k=X_k+e_k+\frac{\mu_k}{L^2},\qquad
 I_i(h_k)=4\int_{\R^3}|\partial_{v_i}\sqrt{h_k}|^2\dd v,
\]
we have $I_i(h_k)\le C(I_{i,k}+e_k)$, and the three-dimensional anisotropic Sobolev inequality yields
\[
 \int_{\R^3}h_k^{5/3}\dd v
 \le Cn_k^{2/3}\bigl(I_1(h_k)I_2(h_k)I_3(h_k)\bigr)^{1/3}
 \le Cn_k^{2/3}\bigl(I_{1,k}+e_k\bigr)^{1/3}Z_k^{2/3}.
\]
Substituting \eqref{eq:generalized-score-longitudinal} for $I_{1,k}$ and using $n_k\le b_k\le\mu_k\le L^2 Z_k$ to cancel $b_k^{-2/3}$, we obtain for every cell $k$
\begin{equation}\label{eq:local-truncated-mixed-nash}
 \int_{\R^3}h_k^{5/3}\dd v
 \le C\left(LZ_k \cD_k^{1/3}+L^{2/3}\mu_k^{1/3}Z_k^{4/3}\right).
\end{equation}

\medskip
\noindent\emph{Step 3.}
Using $\sum_k\chi_k^2=1$ and $g^2/f\le 4h\le 4F$ on $H$, summing $Z_k$ over all cells gives
\begin{equation}\label{eq:truncated-score-cutoff-sum}
 \sum_k Z_k\le C\left(X+\frac{\mu}{L^2}\right).
\end{equation}
Similarly, expanding $\sigma_k(v)-\sigma_k(w)$ in $\cD_k$ and using $\sum_k\chi_k(w)^2|\nabla\chi_k(v)|^2\le CL^{-2}$ together with $|K(v-w)|\le (8\pi|v-w|)^{-1}$, we find
\begin{equation}\label{eq:truncated-dissipation-cutoff-sum}
 \sum_k \cD_k
 \le C\cD+\frac{C}{L^2}\int_{\R^3}\frac{g^2}{f}a[F]\dd v
 \le C\cD+\frac{C}{L^2}\int_{\R^3}h\,a[F]\dd v.
\end{equation}

To bound the potential integral on the right-hand side of \eqref{eq:truncated-dissipation-cutoff-sum} using only the two transverse derivatives in $X$, we project $h$ onto the transverse plane by setting $\bar h(y)=\int_{\R}h(v_1,y)\dd v_1$ for $y=(v_2,v_3)\in\R^2$. Writing $n=\int_{\R^3}h\dd v\le\mu$, Cauchy--Schwarz in $v_1$ and the two-dimensional Gagliardo--Nirenberg inequality for $\sqrt{\bar h}$ yield
\[
 \int_{\R^2}\bar h(y)\dd y=n\le\mu,\qquad
 I_{\R^2}(\bar h)\le X,\qquad
 \|\bar h\|_{L^3(\R^2)}\le Cn^{1/3}X^{2/3}.
\]
Since $|v-w|^{-1}\le |y-z|^{-1}$ for $v=(v_1,y)$ and $w=(w_1,z)$, splitting the two-dimensional Riesz potential at radius $r_0=(n/X)^{1/2}$ gives
\[
 \|a[h]\|_{L^\infty(\R^3)}
 \le \frac{1}{4\pi}\sup_{y\in\R^2}\int_{\R^2}\frac{\bar h(z)}{|y-z|}\dd z
 \le C\Bigl(\|\bar h\|_{L^3(\R^2)}r_0^{1/3}+nr_0^{-1}\Bigr)
 \le C\sqrt{nX}.
\]
By symmetry of the Newton kernel, we obtain
\begin{equation}\label{eq:truncated-cross-potential}
 \int_{\R^3}h\,a[F]\dd v=\int_{\R^3}F\,a[h]\dd v\le C\mu\sqrt{nX}\le C\mu^{3/2}X^{1/2},
\end{equation}
and therefore $\sum_k \cD_k\le C(\cD+L^{-2}\mu^{3/2}X^{1/2})$. Finally, since $h=\sum_k h_k$ with overlap at most $N$, we have $h^{5/3}\le N^{2/3}\sum_k h_k^{5/3}$. Summing \eqref{eq:local-truncated-mixed-nash} over $k$ via H\"older's inequality yields
\[
 \int_{\R^3}h^{5/3}\dd v
 \le C\left[L\Bigl(\sum_k Z_k\Bigr)\Bigl(\sum_k \cD_k\Bigr)^{1/3}
   +L^{2/3}\mu^{1/3}\Bigl(\sum_k Z_k\Bigr)^{4/3}\right],
\]
and optimizing the cell scale with $L=(\mu/X)^{1/2}$ establishes \eqref{eq:truncated-mixed-nash}.
\end{proof}

\begin{remark}
Letting $K\downarrow 0$ in \eqref{eq:truncated-mixed-nash} by monotone convergence yields the global mixed Nash inequality
\[
 \int_{\R^3}f^{5/3}\dd v\le C\left(m^{2/3}I_\perp(f)+m^{1/2}I_\perp(f)^{1/2}\cD(f)^{1/3}\right),
 \qquad m=\int_{\R^3}f\dd v,
\]
where $\cD(f)$ is the full Coulomb entropy dissipation. For the exterior De Giorgi iteration in Section~\ref{sec:original-tail-DG}, the truncated form \eqref{eq:truncated-mixed-nash} is crucial because every quantity on its right-hand side is supported on the upper level set $\{f>K\}$.
\end{remark}

\section{Weighted and localized truncated mixed Nash inequalities}\label{sec:weighted-truncated-nash}

We next extend Proposition~\ref{prop:mixed-nash} to general weak score triplets and localize it to radially aligned velocity cells with polynomial weights. Constants in this section may depend on the upper bound $E_0$ for $\int_{\R^3}\la{v}^2f\dd v$, but are independent of the density threshold $K$ and the tail radius.

\subsection{Formulation for general weak scores}

\begin{proposition}\label{prop:generalized-nash}
Suppose nonnegative integrable functions $F,g,h$ on $\R^3$ satisfy $g\in W^{1,1}_{\mathrm{loc}}(\R^3)$, $\sqrt{h}\in H^1_{\mathrm{loc}}(\R^3)$, and, for a measurable vector field $\sigma$ with $F\sigma\in L^1_{\mathrm{loc}}(\R^3)$,
\[
 h\le g\le F,\qquad \frac{g^2}{F}\le C_0 h,\qquad F\sigma=\nabla g\quad\text{weakly on }\R^3,
\]
where quotients vanish on $\{F=0\}$. Set
\[
 \mu=\int_{\R^3}F\dd v,\qquad X=\int_{\R^3}F|\sigma_\perp|^2\dd v,\qquad
 e=4\int_{\R^3}\left|\nabla\sqrt{h}-\frac{1}{2}\sqrt{F}\,\sigma\right|^2\dd v,\qquad Y=X+e,
\]
and
\[
 \cD=\frac{1}{2}\iint_{\R^3\times\R^3}F(v)F(w)\bigl(\sigma(v)-\sigma(w)\bigr)^\top K(v-w)\bigl(\sigma(v)-\sigma(w)\bigr)\dd v\dd w.
\]
Then there exists a constant $C=C(C_0)>0$ such that
\begin{equation}\label{eq:generalized-mixed-nash}
 \int_{\R^3}h^{5/3}\dd v
 \le C\left(\mu^{2/3}Y+\mu^{1/2}Y^{1/2}\cD^{1/3}\right).
\end{equation}
\end{proposition}

\begin{proof}
We may assume $\mu>0$ and $0<Y+\cD<\infty$, since $I_\perp(h)\le 2(X+e)=2Y$ forces $h=g=0$ whenever $Y=0$. Choose a smooth quadratic partition of unity $\sum_j\chi_j^2=1$ at scale $L>0$ as in Step~2 of Proposition~\ref{prop:mixed-nash}, and set
\[
 F_j=\chi_j^2 F,\qquad g_j=\chi_j^2 g,\qquad h_j=\chi_j^2 h,\qquad
 \sigma_j=\sigma+2\frac{g}{F}\nabla\log\chi_j
\]
on $\{F_j>0\}$, with $\sigma_j=0$ on $\{F_j=0\}$. Then $F_j\sigma_j=\nabla g_j$ weakly on $\R^3$. Writing $\zeta=\nabla\sqrt{h}-\frac{1}{2}\sqrt{F}\,\sigma$, the localized score defect satisfies
\[
 \nabla\sqrt{h_j}-\frac{1}{2}\sqrt{F_j}\,\sigma_j
 =\chi_j\zeta+\left(\sqrt{h}-\frac{g}{\sqrt{F}}\right)\nabla\chi_j,
 \qquad \left(\sqrt{h}-\frac{g}{\sqrt{F}}\right)^2\le (1+C_0)h.
\]
Summing the localized transverse energy $X_j$ and defect $e_j$ gives $\sum_j(X_j+e_j)\le C(Y+\mu/L^2)$. Moreover, since $I_\perp(h)\le 2Y$, the two-dimensional marginal potential estimate \eqref{eq:truncated-cross-potential} gives $\|a[h]\|_{L^\infty(\R^3)}\le C\sqrt{\mu Y}$. Setting $\mu_j=\int_{\R^3} F_j\dd v$ and $Z_j=X_j+e_j+\mu_j/L^2$, the cutoff summation bounds from Step~3 of Proposition~\ref{prop:mixed-nash} yield
\[
 \sum_j \cD_j\le C\left(\cD+L^{-2}\mu^{3/2}Y^{1/2}\right),\qquad
 \sum_j Z_j\le C\left(Y+\frac{\mu}{L^2}\right).
\]
Applying the longitudinal estimate \eqref{eq:generalized-score-longitudinal} and the anisotropic Sobolev inequality on each cell with $\int_{\R^3}g_j\dd v>0$ gives
\[
 \int_{\R^3}h_j^{5/3}\dd v
 \le C\left(LZ_j \cD_j^{1/3}+L^{2/3}\mu_j^{1/3}Z_j^{4/3}\right).
\]
Summing over $j$ and choosing $L=(\mu/Y)^{1/2}$ proves \eqref{eq:generalized-mixed-nash}.
\end{proof}

\subsection{Radially aligned weighted tail estimate}

For $v\in\R^3\setminus\{0\}$, let $\Pi_v=(v\otimes v)/|v|^2$ and $\Pi_v^\perp=I_3-\Pi_v$ denote the radial and transverse projections. The anisotropic background metric $\la{v}^{-3}\Pi_v+\la{v}^{-1}\Pi_v^\perp$ extends smoothly to $\R^3$ as
\[
 \cG(v)=\la{v}^{-1}I_3-\la{v}^{-3}v\otimes v,
 \qquad \cG(0)=I_3.
\]
On a unit-scale velocity cell centered at $v_j$ far from the origin, choosing the local longitudinal axis along the radial unit vector $\omega_j=v_j/|v_j|$ ensures that the two transverse directions lie in the range of $\Pi_{\omega_j}^\perp$, where $\cG(v)$ provides the stronger weight $\la{v_j}^{-1}$. Multiplying each cell inequality by $W(v)=\la{v}^\alpha$ with $\alpha\ge 3/2$ compensates this $\la{v_j}^{-1}$ factor so that the cell Sobolev embedding holds with constants independent of $|v_j|$.

\begin{proposition}\label{prop:weighted-nash}
Let $f\ge 0$ be smooth and rapidly decaying with $\int_{\R^3}\la{v}^2f\dd v\le E_0$ and $A_0=\|a[f]\|_{L^\infty(\R^3)}<\infty$. For $3/2\le \alpha\le 2$, set $W(v)=\la{v}^\alpha$, $P=\sqrt{W}$, $s(v)=\nabla\log f(v)$, and for $K>0$ put
\[
 F=f\mathbf{1}_{\{f>K\}},\qquad g=(f-K)_+,\qquad h=(\sqrt{f}-\sqrt{K})_+^2,
\]
\[
 \cG(v)=\la{v}^{-1}I_3-\la{v}^{-3}v\otimes v,\qquad
 M=\int_{\R^3}WF\dd v,\qquad
 Z=\int_{\R^3}WF\,s^\top\cG s\dd v,
\]
and
\[
 \cD_P=\frac{1}{2}\iint_{\R^3\times\R^3}F(v)F(w)\bigl(P(v)s(v)-P(w)s(w)\bigr)^\top K(v-w)\bigl(P(v)s(v)-P(w)s(w)\bigr)\dd v\dd w.
\]
Then
\begin{equation}\label{eq:weighted-truncated-nash}
 \int_{\R^3}Wh^{5/3}\dd v
 \le C_{E_0}\left[M^{2/3}(Z+M)+M^{1/2}(Z+M)^{1/2}\bigl(\cD_P+Z+(1+A_0)M\bigr)^{1/3}\right].
\end{equation}

More generally, let $\mathcal U\subset\R^3$ be a measurable target set and $\mathcal V$ an open neighborhood containing the closed $c_*\delta$-neighborhood of $\mathcal U$ with $0<\delta\le 1$, where $c_*$ is a fixed numerical constant. Let $p$ be any nonnegative measurable function equal to $P$ on $\mathcal V$, define $\cD_p$ by replacing $P$ with $p$ in $\cD_P$, and set
\[
 M_{\mathcal V}=\int_{\mathcal V}WF\dd v,\qquad
 Z_{\mathcal V}=\int_{\mathcal V}WF\,s^\top\cG s\dd v.
\]
Then
\begin{equation}\label{eq:tail-weighted-truncated-nash}
 \begin{split}
 \int_{\mathcal U}Wh^{5/3}\dd v
 &\le C_{E_0}\Bigl[M_{\mathcal V}^{2/3}\bigl(Z_{\mathcal V}+\delta^{-2}M_{\mathcal V}\bigr)\\
 &\qquad +M_{\mathcal V}^{1/2}\bigl(Z_{\mathcal V}+\delta^{-2}M_{\mathcal V}\bigr)^{1/2}\bigl(\cD_p+Z_{\mathcal V}+\delta^{-2}A_0 M_{\mathcal V}\bigr)^{1/3}\Bigr].
 \end{split}
\end{equation}
In particular, with $B=Z_{\mathcal V}+\cD_p$, the right-hand side of \eqref{eq:tail-weighted-truncated-nash} is bounded by
\[
 C_{E_0,A_0}\delta^{-2}\left[M_{\mathcal V}^{2/3}(B+M_{\mathcal V})+M_{\mathcal V}^{1/2}(B+M_{\mathcal V})^{5/6}\right].
\]
\end{proposition}

\begin{proof}
We may assume $M_{\mathcal V}>0$ and $Z_{\mathcal V}+\cD_p<\infty$. Take a smooth quadratic partition of unity $\{\chi_j\}$ at scale $\delta\in(0,1]$ with uniformly bounded overlap, and retain only those cells whose supports intersect $\mathcal U$. Choosing the numerical constant $c_*>0$ sufficiently large ensures that the support of every retained cell lies inside $\mathcal V\subset\{p=P\}$. Let $v_j$ denote the center of cell $j$. Setting $R_j=\la{v_j}$ and $W_j=R_j^\alpha$, we have on the support of cell $j$
\[
 \la{v}\sim R_j,\qquad W(v)\sim W_j.
\]

For cells with $|v_j|\ge 4$, we align the local longitudinal axis with the radial unit vector $\omega_j=v_j/|v_j|$. Since $|v-v_j|\le C\delta\le C$ on the cell, the angle between $v/|v|$ and $\omega_j$ is at most $C\delta/R_j$, and hence the transverse projection $\Pi_{\omega_j}^\perp=I_3-\omega_j\otimes\omega_j$ satisfies
\[
 R_j^{-1}|\Pi_{\omega_j}^\perp s(v)|^2\le C\,s(v)^\top\cG(v)s(v).
\]
For the remaining cells with $|v_j|<4$, the metric satisfies $\cG(v)\ge cI_3$ throughout the cell support, so the same inequality holds for any unit vector $\omega_j$.

On each retained cell $j$, define the localized triplet $(F_j,g_j,h_j)$ and score $\sigma_j$ as in \eqref{eq:truncated-generalized-score}. Let $\mu_j=\int_{\R^3} F_j\dd v$, and let $X_j$, $e_j$, $Y_j=X_j+e_j$, and $\cD_j$ be the transverse score energy (relative to $\omega_j$), score defect, and cell dissipation of Proposition~\ref{prop:generalized-nash}. Applying Proposition~\ref{prop:generalized-nash} on cell $j$ and setting
\[
 M_j=W_j\mu_j,\qquad Z_j=W_j R_j^{-1}Y_j,\qquad E_j=W_j \cD_j,
\]
we observe that multiplying \eqref{eq:generalized-mixed-nash} by $W_j=R_j^\alpha$ produces the radial scaling powers $R_j^{1-2\alpha/3}$ and $R_j^{1/2-\alpha/3}$ on the right-hand side, both of which are bounded above by $1$ because $\alpha\ge 3/2$. Consequently,
\[
 W_j\int_{\R^3}h_j^{5/3}\dd v
 \le C\left(M_j^{2/3}Z_j+M_j^{1/2}Z_j^{1/2}E_j^{1/3}\right).
\]
Using \eqref{eq:truncated-gradient-error} and the transverse metric bound above, summing the single-particle quantities over all retained cells yields
\begin{equation}\label{eq:weighted-cell-Fisher-sum}
 \sum_j M_j\le CM_{\mathcal V},\qquad
 \sum_j Z_j\le C\bigl(Z_{\mathcal V}+\delta^{-2}M_{\mathcal V}\bigr).
\end{equation}

To bound the weighted cell dissipation sum $\sum_j E_j$, note that any pair $(v,w)$ in the support of cell $j$ lies in $\mathcal V\times\mathcal V$ and satisfies $|v-w|\le C\delta$. The cutoff gradient contribution to $\sum_j E_j$ is bounded by
\[
 C\delta^{-2}\int_{\mathcal V}Wh\,a[F]\dd v\le C\delta^{-2}A_0 M_{\mathcal V}.
\]
For the unweighted score difference $s(v)-s(w)$ on cell $j$, we use $p=P$ on $\mathcal V$ to write
\[
 P(v)\bigl(s(v)-s(w)\bigr)=\bigl(p(v)s(v)-p(w)s(w)\bigr)-\bigl(P(v)-P(w)\bigr)s(w).
\]
The first term on the right-hand side sums to at most $C\cD_p$. For the commutator term involving $P(v)-P(w)$, the mean value theorem and Chebyshev's inequality with the second moment give
\[
 |P(v)-P(w)|\le CR_j^{\alpha/2-1}|v-w|,\qquad
 \int_{\{|v-w|\le C\delta\}}F(v)\dd v\le CE_0\la{w}^{-2}.
\]
Summing over $j$ and using $|K(v-w)|\le (8\pi|v-w|)^{-1}$ therefore bounds the commutator contribution by
\begin{align*}
 &C\int_{\mathcal V}F(w)|s(w)|^2\la{w}^{\alpha-2}\left(\int_{\{|v-w|\le C\delta\}}F(v)|v-w|\dd v\right)\dd w\\
 &\qquad \le CE_0\delta\int_{\mathcal V}F(w)\la{w}^{\alpha-4}|s(w)|^2\dd w
 \le CE_0 Z_{\mathcal V},
\end{align*}
where the last inequality uses $\cG(w)\ge\la{w}^{-3}I_3$ and $W(w)=\la{w}^\alpha$. Combining these estimates gives
\begin{equation}\label{eq:weighted-cell-dissipation-sum}
 \sum_j E_j\le C_{E_0}\bigl(\cD_p+Z_{\mathcal V}+\delta^{-2}A_0 M_{\mathcal V}\bigr).
\end{equation}
Finally, summing the cell Nash inequalities over $j$ via H\"older's inequality and substituting \eqref{eq:weighted-cell-Fisher-sum}--\eqref{eq:weighted-cell-dissipation-sum} establishes \eqref{eq:tail-weighted-truncated-nash}. Taking $\mathcal U=\mathcal V=\R^3$, $p=P$, and $\delta=1$ gives \eqref{eq:weighted-truncated-nash}.
\end{proof}

\section{Weighted tail entropy and proof of Theorem~\ref{thm:global-bound}}\label{sec:original-tail-DG}

Let $f$ be a smooth, rapidly decaying solution of \eqref{eq:landau}--\eqref{eq:normalization} on $[0,T]\times\R^3$, fix $0<\tau<T$, and write $s=\nabla\log f$, $A=A[f]$, and $b=\nabla a[f]$. By Theorem~\ref{thm:critical-fisher}, Corollary~\ref{cor:smoothing}, and \eqref{eq:weighted-L3}, we have the uniform bounds
\begin{equation}\label{eq:tail-preparatory-bounds}
 \sup_{\tau/4\le t\le T}
 \left\{I_{-1}(f(t))+\|a[f(t)]\|_{L^\infty(\R^3)}
       +\int_{\R^3}\la{v}^{7/4}f(t,v)^{11/10}\dd v\right\}
 \le C_{H,\tau},
\end{equation}
where the weighted $L^{11/10}$ estimate follows by interpolating $\|\la{v}^{-1}f\|_{L^3}$ with the conserved energy:
\[
 \int_{\R^3}\la{v}^{7/4}f^{11/10}\dd v
 \le \left(\int_{\R^3}\la{v}^{-3}f^3\dd v\right)^{1/20}
      \left(\int_{\R^3}\la{v}^2f\dd v\right)^{19/20}.
\]
Furthermore, setting $\tilde f(w)=\la{w}^{-1}f(w)$ and splitting the convolution for $b(t,v)=\nabla a[f](t,v)$ at $|v-w|=\la{v}/2$, the translated Hardy inequality \eqref{eq:Hardy} on $\{|v-w|\le\la{v}/2\}$ gives
\[
 \int_{\{|v-w|\le \la{v}/2\}}\frac{f(w)}{|v-w|^2}\dd w
 \le C\la{v}\int_{\R^3}\frac{\tilde f(w)}{|v-w|^2}\dd w
 \le C\la{v} I(\tilde f)
 \le C\la{v}\bigl(1+I_{-1}(f)\bigr),
\]
while the complementary integral over $\{|v-w|>\la{v}/2\}$ is bounded by $4\la{v}^{-2}$. Consequently, the drift field satisfies the linear growth bound
\begin{equation}\label{eq:linear-drift-tail}
 |b(t,v)|\le C_{H,\tau}\la{v}
 \qquad \bigl((t,v)\in[\tau/4,T]\times\R^3\bigr).
\end{equation}

We also record the anisotropic lower ellipticity of the diffusion matrix generated by the low-density part of $f$. There exists a threshold $K_0=K_0(H)\ge 2$ depending only on $H$ such that, for every $k\ge K_0$ and all $(t,v)\in[0,T]\times\R^3$,
\begin{equation}\label{eq:low-density-anisotropy-tail}
 A\bigl[f\mathbf{1}_{\{f\le k\}}\bigr](t,v)
 \ge c_H\bigl(\la{v}^{-3}\Pi_v+\la{v}^{-1}\Pi_v^\perp\bigr)
 =c_H\cG(v).
\end{equation}
Indeed, using $\int_{\R^3}f\log_+ f\dd w\le C_H$ and $\int_{\R^3}|w|^2f(w)\dd w=3$, we choose $K_0\ge 2$ large enough that $C_H/\log K_0\le 1/4$. Then the truncated density
\[
 q(w)=f(w)\mathbf{1}_{\{|w|\le4\}}\mathbf{1}_{\{f(w)\le K_0\}}
\]
has mass $m_q\ge 1-3/16-1/4=9/16$, is supported in $B_4$, and satisfies $q\le K_0$, so its center of mass $\bar w=m_q^{-1}\int_{\R^3}wq(w)\dd w$ lies in $\overline B_4$. For any unit vector $e\in S^2$, the slab $\{w\in B_4:|e\cdot(w-\bar w)|<d_H\}$ has measure at most $32\pi d_H$, and choosing $d_H=(32\pi K_0)^{-1}(9/32)$ ensures that at least $9/32$ of the mass of $q$ lies outside this slab, whence
\[
 \int_{\R^3} q(w)\bigl(e\cdot(w-\bar w)\bigr)^2\dd w\ge \frac{9}{32}d_H^2.
\]
Applying this lower bound to two orthogonal unit vectors perpendicular to $\xi\in S^2$ and expanding about $\bar w$ yields
\[
 \int_{\R^3} q(w)|\xi\times(v-w)|^2\dd w
 =m_q|\xi\times(v-\bar w)|^2
   +\int_{\R^3} q(w)|\xi\times(w-\bar w)|^2\dd w
 \ge c_H\bigl(|\xi\times v|^2+1\bigr).
\]
Since $|v-w|\le C\la{v}$ on $B_4$ and $q\le f\mathbf{1}_{\{f\le k\}}$ for all $k\ge K_0$, dividing by $8\pi|v-w|^3$ establishes \eqref{eq:low-density-anisotropy-tail}.

\subsection{An exact weighted truncated entropy identity}

For each threshold $k>0$, we define the upper level set truncations and the truncated relative entropy density $\psi_k(f)$ by
\[
 H_k=\{f>k\},\qquad F_k=f\mathbf{1}_{H_k},\qquad
 g_k=(f-k)_+,\qquad h_k=(\sqrt f-\sqrt k)_+^2,
\]
and
\[
 \psi_k(f)=\bigl(f\log(f/k)-f+k\bigr)\mathbf{1}_{H_k}.
\]
Given a smooth nonnegative spatial weight $p(v)$, we define the weighted truncated entropy $E_{p,k}$ and the split dissipation functional $B_{p,k}=\cD_{p,k}+L_{p,k}$ by
\[
 E_{p,k}=\int_{\R^3}p^2\psi_k(f)\dd v,\qquad
 L_{p,k}=\int_{H_k}p^2f\,s^\top A\bigl[f\mathbf{1}_{\{f\le k\}}\bigr]s\dd v,
\]
and
\[
 \cD_{p,k}=\frac12\iint_{\R^3\times\R^3} F_k(v)F_k(w)
 \bigl(p(v)s(v)-p(w)s(w)\bigr)^\top K(v-w)
 \bigl(p(v)s(v)-p(w)s(w)\bigr)\dd v\dd w.
\]

\begin{lemma}\label{lem:weighted-entropy}
Writing $b_K(z)=\diver K(z)=-z/(4\pi|z|^3)$, the weighted truncated entropy satisfies the exact evolution identity
\begin{equation}\label{eq:exact-weighted-tail-entropy}
 \frac{\dd}{\dd t}E_{p,k}+B_{p,k}
 =k\int_{\R^3} p^2g_k\dd v
  +\int_{\R^3}\psi_k(f)\bigl(A:D^2(p^2)+2b\cdot\nabla(p^2)\bigr)\dd v
  -\cE_p(g_k),
\end{equation}
where the nonlocal weight commutator $\cE_p(g)$ is given by
\begin{equation}\label{eq:spatial-weight-cross-error}
 \begin{split}
 \cE_p(g)=\iint_{\R^3\times\R^3} g(v)g(w)\Bigl[&\nabla p(v)^\top K(v-w)\nabla p(w)\\
 &+b_K(v-w)\cdot\bigl(p(v)\nabla p(w)-p(w)\nabla p(v)\bigr)\Bigr]\dd v\dd w.
 \end{split}
\end{equation}
\end{lemma}

\begin{proof}
The truncated entropy density $\psi_k$ belongs to $C^1([0,\infty))$ and satisfies
\[
 \psi_k'(x)=\log(x/k)\mathbf{1}_{\{x>k\}},\qquad
 \psi_k''(x)=x^{-1}\mathbf{1}_{\{x>k\}}\quad (x\ne k),\qquad
 x^2\psi_k'(x)-x\psi_k(x)=xg_k(x).
\]
Differentiating $E_{p,k}$ along the non-divergence equation $\partial_t f=A:D^2f+f^2$, integrating by parts twice in $v$, and using the identity $\diver\diver(p^2A)=A:D^2(p^2)+2b\cdot\nabla(p^2)-p^2f$ yield
\[
 \frac{\dd}{\dd t}E_{p,k}
 =-\int_{H_k}p^2f\,s^\top As\dd v+\int_{\R^3} p^2fg_k\dd v
   +\int_{\R^3}\psi_k(f)\bigl(A:D^2(p^2)+2b\cdot\nabla(p^2)\bigr)\dd v.
\]
Splitting the diffusion matrix $A[f]=A[f\mathbf{1}_{\{f\le k\}}]+A[F_k]$ separates the low-density diffusion $L_{p,k}$ from the high-density diffusion integral. Using the weak gradient identity $F_k s=\nabla g_k$ and expanding the square in $\cD_{p,k}$ by symmetry in $(v,w)$, we find
\[
 \int_{H_k}p^2f\,s^\top A[F_k]s\dd v=\cD_{p,k}+J_{p,k},
\]
where the high-density cross term is
\[
 J_{p,k}=\iint_{\R^3\times\R^3} p(v)p(w)\nabla g_k(v)^\top K(v-w)\nabla g_k(w)\dd v\dd w.
\]
In the sense of distributions on $\R^3$, the Landau kernel satisfies
\[
 K(z)=\frac{1}{8\pi}D^2|z|,\qquad
 \diver K(z)=b_K(z),\qquad \diver b_K=-\delta_0.
\]
Integrating $J_{p,k}$ by parts in both $v_i$ and $w_j$ and using $\partial_{w_j}K_{ij}(v-w)=-\partial_{z_j}K_{ij}(z)$ give
\[
 \begin{split}
 \partial_{v_i}\partial_{w_j}\bigl(p(v)p(w)K_{ij}(v-w)\bigr)
 &=\nabla p(v)^\top K(v-w)\nabla p(w)\\
 &\quad+b_K(v-w)\cdot\bigl(p(v)\nabla p(w)-p(w)\nabla p(v)\bigr)
       +p(v)p(w)\delta_0(v-w),
 \end{split}
\]
and therefore
\[
 J_{p,k}=\int_{\R^3} p^2g_k^2\dd v+\cE_p(g_k).
\]
On compact sets, $|p(v)\nabla p(w)-p(w)\nabla p(v)|=O(|v-w|)$ cancels one power of the $|v-w|^{-2}$ singularity in $b_K(v-w)$, so the kernel of $\cE_p(g_k)$ is locally integrable and the distributional integration by parts is justified by standard mollification of the compactly supported Lipschitz function $g_k$. Finally, substituting $J_{p,k}$ into the evolution equation and using $fg_k-g_k^2=kg_k$ establishes \eqref{eq:exact-weighted-tail-entropy}.
\end{proof}

\subsection{Nested exterior cutoffs and initial tail smallness}

Fix a base radius $R\ge 2$ to be chosen below, and define for each integer $n\ge 0$ the nested radii, density thresholds, and time levels
\[
 R_n=R(2-2^{-n}),\qquad k_n=K_0(2-2^{-n}),\qquad
 t_n=\frac{\tau}{2}+\frac{\tau}{4}(1-2^{-n}).
\]
Choose smooth radial cutoffs $0\le\chi_n(v)\le 1$ satisfying
\[
 \chi_n(v)=0\quad (|v|\le R_n),\qquad
 \chi_n(v)=1\quad (|v|\ge R_{n+1}),\qquad
 |\nabla\chi_n|\le\frac{C2^n}{R},\qquad
 |D^2\chi_n|\le\frac{C4^n}{R^2},
\]
and set
\[
 W(v)=\la{v}^{3/2},\quad P(v)=\la{v}^{3/4},\quad p_n(v)=P(v)\chi_n(v),
 \qquad E_n=E_{p_n,k_n},\quad B_n=B_{p_n,k_n}.
\]
Throughout this section, $b_0>1$ denotes a sufficiently large universal geometric base (such as $b_0=16$) independent of $R$, $n$, and $f$. Since $\la{v}\sim R$ on the cutoff transition annulus $\{R_n\le |v|\le R_{n+1}\}$ and $R\ge 2$, the product rule yields the derivative bounds
\[
 \begin{aligned}
 |\nabla p_n|&\le Cb_0^n\la{v}^{-1/4}\mathbf{1}_{\{|v|>R_n\}},&
 |D^2p_n|&\le Cb_0^n\la{v}^{-5/4}\mathbf{1}_{\{|v|>R_n\}},\\
 |\nabla(p_n^2)|&\le Cb_0^n\la{v}^{1/2}\mathbf{1}_{\{|v|>R_n\}},&
 |D^2(p_n^2)|&\le Cb_0^n\la{v}^{-1/2}\mathbf{1}_{\{|v|>R_n\}}.
 \end{aligned}
\]
Combining these bounds with $\|A[f]\|\le\tr A[f]=a[f]$, \eqref{eq:tail-preparatory-bounds}, and the linear drift estimate \eqref{eq:linear-drift-tail} gives
\begin{equation}\label{eq:weight-generator-tail-bound}
 \bigl|A:D^2(p_n^2)+2b\cdot\nabla(p_n^2)\bigr|
 \le C_{H,\tau}b_0^n W\mathbf{1}_{\{|v|>R_n\}}.
\end{equation}

We next show that for any $0\le g\le f$, the nonlocal weight commutator $\cE_{p_n}(g)$ in \eqref{eq:spatial-weight-cross-error} satisfies
\begin{equation}\label{eq:tail-cross-error-bound}
 |\cE_{p_n}(g)|
 \le Cb_0^n\bigl(R^{-2}\|a[f]\|_{L^\infty(\R^3)}+R^{-3}\bigr)
           \int_{\{|v|>R_n\}}Wg\dd v.
\end{equation}
To verify \eqref{eq:tail-cross-error-bound}, let $\mathcal K_n(v,w)$ denote the bracketed kernel in \eqref{eq:spatial-weight-cross-error}, which is supported on $\{|v|>R_n\}\times\{|w|>R_n\}$. For nearby pairs with $|v-w|\le\frac{1}{8}\max(\la{v},\la{w})$, all points on the line segment between $v$ and $w$ have velocities comparable to $\la{v}\sim\la{w}$, so the mean value theorem and the derivative bounds for $p_n$ give
\[
 |\nabla p_n(v)|\,|\nabla p_n(w)|\le Cb_0^n\la{v}^{-1/2},\qquad
 |p_n(v)\nabla p_n(w)-p_n(w)\nabla p_n(v)|\le Cb_0^n\la{v}^{-1/2}|v-w|.
\]
Thus $|\mathcal K_n(v,w)|\le Cb_0^n\la{v}^{-1/2}|v-w|^{-1}$, and since $\la{v}^{-1/2}\le R^{-2}W(v)$ on $\{|v|>R_n\}$, the integral over nearby pairs is bounded by
\[
 Cb_0^n\int_{\{|v|>R_n\}}\la{v}^{-1/2}g(v)
           \left(\int_{\R^3}\frac{g(w)}{|v-w|}\dd w\right)\dd v
 \le Cb_0^nR^{-2}\|a[f]\|_{L^\infty(\R^3)}\int_{\{|v|>R_n\}}Wg\dd v.
\]
For separated pairs with $|v-w|>\frac{1}{8}\max(\la{v},\la{w})$, symmetry allows us to assume $\la{v}\ge\la{w}>R$. The bounds on $p_n$ and $\nabla p_n$ then give
\[
 |\mathcal K_n(v,w)|
 \le Cb_0^n\la{v}^{-5/4}\la{w}^{-1/4}\le Cb_0^nR^{-3/2},
\]
and using $\int_{\{|v|>R_n\}}g\dd v\le\min(1,R^{-3/2}\int_{\{|v|>R_n\}}Wg\dd v)$ bounds the integral over separated pairs by $Cb_0^nR^{-3}\int_{\{|v|>R_n\}}Wg\dd v$, proving \eqref{eq:tail-cross-error-bound}.

Substituting \eqref{eq:weight-generator-tail-bound} and \eqref{eq:tail-cross-error-bound} into Lemma~\ref{lem:weighted-entropy} with $p=p_n$ and $k=k_n\le 2K_0$ gives, on $[\tau/4,T]$,
\begin{equation}\label{eq:tail-energy-inequality}
 \frac{\dd}{\dd t}E_n+B_n
 \le C_{H,\tau}b_0^n\int_{\{|v|>R_n\}}
             W\bigl(\psi_{k_n}(f)+f\mathbf{1}_{\{f>k_n\}}\bigr)\dd v.
\end{equation}
For $k\in[K_0,2K_0]$, the scalar inequality $\psi_k(x)\le C_{K_0}x^{11/10}$ together with \eqref{eq:tail-preparatory-bounds} and the conserved energy $\int_{\R^3}\la{v}^2f\dd v=4$ yields the uniform tail decay estimates
\begin{equation}\label{eq:tail-static-smallness}
 \sup_{\tau/4\le t\le T}\int_{\{|v|>R\}}W\psi_k(f)\dd v
 \le C_{H,\tau}R^{-1/4},\qquad
 \sup_{0\le t\le T}\int_{\{|v|>R\}}Wf\dd v\le CR^{-1/2}.
\end{equation}
Defining the space-time tail energy functionals
\begin{equation}\label{eq:tail-initial-energies}
 Y_n=\sup_{t_n\le t\le T}E_n(t)+\int_{t_n}^TB_n(t)\dd t,
 \qquad Y_0,Y_1\le C_{H,\tau,T}R^{-1/4},
\end{equation}
the initial smallness bounds for $Y_0$ and $Y_1$ follow immediately by integrating \eqref{eq:tail-energy-inequality} from $t=\tau/4$ and applying \eqref{eq:tail-static-smallness}.

\subsection{De Giorgi recurrence and completion of the proof}

For each step $n\ge 2$, define the weighted mass and $L^{5/3}$ truncated norm
\[
 \mu_n(t)=\int_{\{|v|>R_{n-1}\}}Wf\mathbf{1}_{\{f>k_{n-1}\}}\dd v,
 \qquad
 \Phi_n(t)=\int_{\{|v|>R_n\}}Wh_{k_{n-1}}^{5/3}\dd v.
\]
We apply Proposition~\ref{prop:weighted-nash} at threshold $k_{n-1}$ with exponent $\alpha=3/2$, target set $\mathcal U=\{|v|>R_n\}$, neighborhood $\mathcal V=\{|v|>R_{n-1}\}$, weight $p=p_{n-2}$ (which equals $P=W^{1/2}$ on $\mathcal V$), and cell scale
\[
 \delta_n=\min\bigl(1,(2c_*)^{-1}R2^{-n}\bigr),\qquad
 \delta_n^{-2}\le C4^n.
\]
To bound the dissipation terms at threshold $k_{n-1}$ by $B_{n-2}$ at threshold $k_{n-2}$, we compare the two components of $B_{p,k}$ separately: at the fixed weight $p=p_{n-2}$, the high-density dissipation satisfies $\cD_{p_{n-2},k_{n-1}}\le\cD_{p_{n-2},k_{n-2}}$ because $F_{k_{n-1}}\le F_{k_{n-2}}$, while the low-density ellipticity bound \eqref{eq:low-density-anisotropy-tail} at level $k_{n-2}$ gives
\[
 \begin{split}
 Z_{\mathcal V,k_{n-1}}
 &=\int_{\mathcal V}WF_{k_{n-1}}\,s^\top\cG s\dd v\\
 &\le c_H^{-1}\int_{\mathcal V}p_{n-2}^2F_{k_{n-1}}
          s^\top A\bigl[f\mathbf{1}_{\{f\le k_{n-2}\}}\bigr]s\dd v
 \le c_H^{-1}L_{p_{n-2},k_{n-2}}.
 \end{split}
\]
Consequently, $Z_{\mathcal V,k_{n-1}}+\cD_{p_{n-2},k_{n-1}}\le C_H B_{n-2}$, and Proposition~\ref{prop:weighted-nash} yields
\begin{equation}\label{eq:nested-static-input-tail}
 \Phi_n\le C_{H,\tau}b_0^n\left[
     \mu_n^{2/3}(B_{n-2}+\mu_n)
    +\mu_n^{1/2}(B_{n-2}+\mu_n)^{5/6}\right].
\end{equation}

On the dyadic thresholds $k_n=K_0(2-2^{-n})\in[K_0,2K_0]$, elementary calculus on $\psi_k(x)$ and $h_k(x)=(\sqrt{x}-\sqrt{k})_+^2$ provides the level comparisons
\begin{equation}\label{eq:tail-level-comparisons}
 \begin{split}
 f\mathbf{1}_{\{f>k_{n-1}\}}&\le C_{K_0}4^n\psi_{k_{n-2}}(f),\\
 \psi_{k_n}(f)+f\mathbf{1}_{\{f>k_n\}}
 &\le C_{K_0}b_0^n h_{k_{n-1}}^{5/3},\qquad n\ge2.
 \end{split}
\end{equation}
Indeed, for $x>k_{n-1}>k_{n-2}$, we have
\[
 \psi_{k_{n-2}}(x)\ge\frac{(x-k_{n-2})^2}{2x}\ge 2^{-2n-1}x,
\]
which gives the first inequality in \eqref{eq:tail-level-comparisons}; for $x>k_n$, the lower bounds
\[
 \sqrt{x}-\sqrt{k_{n-1}}\ge c\sqrt{K_0}\,2^{-n}\quad (k_n<x\le 4K_0),\qquad
 h_{k_{n-1}}(x)\ge cx\quad (x\ge 4K_0)
\]
give the second inequality.

Fix a smooth nondecreasing time cutoff $0\le\eta_n(t)\le 1$ equal to $0$ on $[0,t_{n-1}]$ and $1$ on $[t_n,T]$, with $0\le\eta_n'\le C2^n/\tau$. Multiplying \eqref{eq:tail-energy-inequality} by $\eta_n(t)$, applying the second inequality in \eqref{eq:tail-level-comparisons}, and integrating in time from $t_{n-1}$ yields
\begin{equation}\label{eq:tail-DG-energy-step}
 Y_n\le C_{H,\tau,T}b_0^n\int_{t_{n-1}}^T\Phi_n(t)\dd t.
\end{equation}
On the other hand, since $p_{n-2}^2=W$ on $\{|v|>R_{n-1}\}$ and $t_{n-1}>t_{n-2}$, the first inequality in \eqref{eq:tail-level-comparisons} implies that on $[t_{n-1},T]$,
\[
 \mu_n(t)\le C_{K_0}4^n E_{n-2}(t)\le C_{K_0}4^n Y_{n-2}.
\]
Setting $U_n=\max(Y_{n-1},Y_{n-2})$ and applying H\"older's inequality in time to \eqref{eq:nested-static-input-tail} and \eqref{eq:tail-DG-energy-step}, we arrive at the superlinear recurrence
\begin{equation}\label{eq:tail-final-recurrence}
 Y_n\le C\gamma_0^n\left(U_n^{5/3}+U_n^{4/3}\right),\qquad n\ge2,
\end{equation}
where $C=C(H,\tau,T)>0$ and the universal geometric base $\gamma_0>1$ are independent of $R$.

Choosing $0<\varepsilon\le 1$ small enough that $2C\varepsilon^{1/3}\gamma_0^{16}\le 1$, and then fixing the radius $R=R(H,\tau,T)\ge 2$ in \eqref{eq:tail-initial-energies} large enough that
\[
 Y_0,Y_1\le\varepsilon \gamma_0^{-6},
\]
induction on $n\ge 0$ using \eqref{eq:tail-final-recurrence} shows that $Y_n\le\varepsilon\gamma_0^{-6n}$ for all $n\ge 0$, and hence $\lim_{n\to\infty}Y_n=0$. For every $t\in[3\tau/4,T]$ and every $n\ge 0$, we have $t\ge t_n$ and $p_n^2=W$ on $\{|v|>2R\}$. Passing to the limit $n\to\infty$ via Fatou's lemma gives
\[
 0\le\int_{\{|v|>2R\}}W\psi_{2K_0}(f(t,v))\dd v
 \le\liminf_{n\to\infty}E_n(t)\le\lim_{n\to\infty}Y_n=0,
\]
which implies $f(t,v)\le 2K_0(H)$ for all $(t,v)\in[3\tau/4,T]\times\{|v|\ge 2R\}$. Combining this uniform exterior tail bound with the interior Moser bound \eqref{eq:local-infty} on the ball $B_{2R(H,\tau,T)}$ yields
\[
 \sup_{\tau\le t\le T}\|f(t)\|_{L^\infty(\R^3)}
 \le \max\bigl\{2K_0(H),\,C\bigl(H,\tau,T,2R(H,\tau,T)\bigr)\bigr\}
 =C(H,\tau,T),
\]
completing the proof of Theorem~\ref{thm:global-bound}.

\medskip
\noindent\textbf{Disclosure of AI assistance.}
During the preparation of this manuscript, the author had useful discussions with OpenAI's ChatGPT and Google's Gemini, and used them to assist with language polishing and checking auxiliary calculations. The author takes full responsibility for all proofs, as well as the accuracy, originality, and integrity of the entire paper.


\begin{thebibliography}{99}
\sloppy

\bibitem{Desvillettes}
L. Desvillettes,
\emph{Entropy dissipation estimates for the Landau equation in the Coulomb case and applications},
J. Funct. Anal. \textbf{269} (2015), no.~5, 1359--1403.

\bibitem{DGGL}
L. Desvillettes, W. Golding, M.~P. Gualdani and A. Loher,
\emph{Production of the Fisher information for the Landau--Coulomb equation with $L^1$ initial data},
arXiv preprint arXiv:2410.10765, 2024.

\bibitem{FW}
X. Feng and Z. Wang,
\emph{Kac's program for the Landau equation},
Arch. Ration. Mech. Anal. \textbf{250} (2026), article~72.

\bibitem{Fournier}
N. Fournier,
\emph{Uniqueness of bounded solutions for the homogeneous Landau equation with a Coulomb potential},
Comm. Math. Phys. \textbf{299} (2010), no.~3, 765--782.

\bibitem{GGL}
W. Golding, M. Gualdani and A. Loher,
\emph{Global smooth solutions to the Landau--Coulomb equation in $L^{3/2}$},
Arch. Ration. Mech. Anal. \textbf{249} (2025), article~34.

\bibitem{GIJV}
F. Golse, C. Imbert, S. Ji and A.~F. Vasseur,
\emph{Local regularity for the space-homogeneous Landau equation with very soft potentials},
arXiv preprint arXiv:2206.05155, 2022, revised 2024.

\bibitem{GS}
N. Guillen and L. Silvestre,
\emph{The Landau equation does not blow up},
Acta Math. \textbf{234} (2025), no.~2, 315--375.

\bibitem{JiSphere}
S. Ji,
\emph{Bounds for the optimal constant of the Bakry--\'Emery $\Gamma_2$ criterion inequality on $\mathbb{RP}^{d-1}$},
arXiv preprint arXiv:2408.13954, 2024.

\bibitem{JiDissipation}
S. Ji,
\emph{Dissipation estimates of the Fisher information for the Landau equation},
arXiv preprint arXiv:2410.09035, 2024.

\bibitem{Silvestre}
L. Silvestre,
\emph{Upper bounds for parabolic equations and the Landau equation},
J. Differential Equations \textbf{262} (2017), no.~3, 3034--3055.

\bibitem{Wu}
K.-C. Wu,
\emph{Global in time estimates for the spatially homogeneous Landau equation with soft potentials},
J. Funct. Anal. \textbf{266} (2014), no.~5, 3134--3155.

\end{thebibliography}
\end{document}